\documentstyle[12pt,a4]{article}
\newcommand{\be}{\begin{equation}}
\newcommand{\ee}{\end{equation}}
\newcommand{\ba}{\begin{eqnarray}}
\newcommand{\ea}{\end{eqnarray}}
\newcommand{\baa}{\begin{eqnarray*}}
\newcommand{\eaa}{\end{eqnarray*}}
\newcommand{\bb}{\begin{thebibliography}}
\newcommand{\eb}{\end{thebibliography}}

\newcommand{\bi}[1]{\bibitem{#1}}
\newcommand{\lab}[1]{\label{#1}}
\newcommand{\re}[1]{(\ref{#1})}

\newcounter{my}
\newcommand{\he}%
    {\stepcounter{equation}\setcounter{my}%
    {\value{equation}}\setcounter{equation}0%
    }%
\newcommand{\she}%
    {\setcounter{equation}{\value{my}}%
     }%

\renewcommand\t{\tilde}

\newcommand\cn{\mbox{cn}}
\newcommand\dn{\mbox{dn}}
\newcommand\sn{\mbox{sn}}

\newcommand\dc{\mbox{dc}}

\newtheorem{th}{Theorem}
\newtheorem{pr}{Proposition}

\newtheorem{lem}{Lemma}

\begin{document}

%\begin{titlepage}
\vspace*{10mm}
\begin{center}

{\LARGE \bf Elliptic solutions of the Toda chain and a generalization of the Stieltjes-Carlitz
polynomials}

\par
\vspace{5mm}
{\large \bf Alexei Zhedanov}
\par
\medskip
{\em Donetsk Institute for Physics and Technology, Donetsk 83114,
Ukraine}
\end{center}

\vspace*{5mm}

\begin{abstract}
We construct new elliptic solutions of the restricted Toda chain.
These solutions give rise to a new explicit class of orthogonal
polynomials which can be considered as a generalization of the
Stieltjes-Carlitz elliptic polynomials. Relations between characteristic (i.e. positive definite) functions, Toda chain and orthogonal polynomials are developed in order to derive main properties of these polynomials. The recurrence
coefficients and the weight function of these polynomials are
expressed explicitly. In the degenerated cases of the elliptic
functions the modified Meixner polynomials and the Krall-Laguerre
polynomials appear.

\end{abstract}

%\end{titlepage}
\newpage
\section{Restricted Toda chain and orthogonal polynomials}
\setcounter{equation}{0} The main purpose of the present paper is
explicit construction of orthogonal polynomials which are a
generalization of the famous Stieltjes-Carlitz orthogonal
polynomials connected with elliptic functions. For theory of these
polynomials and history of their discovery see, e.g.
\cite{Rogers}, \cite{LoBr}, \cite{Carlitz}, \cite{Valent}.

Our main tool will be connection between orthogonal polynomials
depending on an additional "time" parameter $t$ and solutions of
the so-called restricted Toda chain. This connection allows us
first to construct some explicit "elliptic" solutions of the Toda
chain and then to reconstruct corresponding orthogonal polynomials
and their orthogonality measure. As we will see, the
Stieltjes-Carlitz polynomials appear to be a very special case of
these constructed orthogonal polynomials corresponding to a "zero
time" $t=0$ case. We derive also explicit recurrence coefficients
$u_n(t), b_n(t)$ for the obtained orthogonal polynomials. These
coefficients are expressed in terms of the Weierstrass elliptic
functions in $n$ and $t$ . Recall that the recurrence coefficients
of the Stieltjes-Carlitz polynomials are expressed in terms of
linear and quadratic polynomials in $n$. We show that obtained
polynomials are orthogonal on the whole real axis with a positive
discrete measure. The measure is constructed explicitly. It
appears that the parameter $t$ can take values only inside of some
interval (so-called admissible interval) in order for the measure
and the recurrence coefficients will be well defined. The
Stieltjes-Carlitz case corresponds to the middlepoint $t=0$ of
this interval.

We recall basic definitions and results concerning relations
between Toda chain and orthogonal polynomials \cite{Sogo}, \cite{Apt}, \cite{NaZhe}.

The Toda chain equations are \cite{Toda} \be \dot u_n = u_n (b_n -
b_{n-1}), \quad \dot b_n = u_{n+1}-u_n \lab{Toda} \ee with
additional condition \be u_0=0\quad \lab{res} \ee

where the dot indicates the differentiation with respect to $t$.
In what follows we will call equations \re{Toda} with restriction
\re{res} {\it the restricted Toda chain}  (TC) equations.
\par
Let $P_n(x;t)$ be orthogonal polynomials satisfying the three-term
recurrence relation \be P_{n+1}(x) + b_n P_n(x) + u_n P_{n-1}(x) =
x P_n(x) \lab{3_term} \ee with initial conditions \be P_0=1, \quad
P_1(x) = x-b_0. \lab{ini} \ee We will assume that $u_n \ne 0,
n=1,2,\dots$. By the well known spectral theorem \cite{Chi}, \cite{Ismail} there exists a
nondegenerate linear functional $\sigma$ such that the polynomials
$P_n(x)$ are orthogonal with respect to it: \be \sigma(P_n(x)
P_m(x) ) = h_n \delta_{nm}, \lab{ort} \ee where $h_n$ are
normalization constants. The linear functional $\sigma$ can be
defined through its moments \be c_n = \sigma(x^n), \quad
n=0,1,\dots. \lab{moms} \ee It is usually assumed that $c_0 =1$
(standard normalization condition), but we will not assume this
condition in the followings. So we will assume that $c_0$ is an
arbitrary nonzero parameter.
\par
Introduce the Hankel determinants \be D_n = \det
(c_{i+j})_{i,j=0,\dots,n-1}, \quad D_0 =1,\quad D_1= c_0.
\lab{Delta} \ee Then the polynomials $P_n(x)$ can be uniquely
represented as \cite{Chi} \ba P_n(x)=\frac{1}{D_n}\left |
\begin{array}{cccc} c_0 & c_{1} & \dots &
c_{n}\\ c_{1}& c_{2} & \dots & c_{n+1}\\ \dots & \dots & \dots & \dots\\
c_{n-1} & c_{n} & \dots & c_{2n-1} \\ 1 & x & \dots & x^n
\end{array} \right |. \lab{P_det} \ea The normalization constants
are expressed as \be h_n = \frac{D_{n+1}}{D_n}, \quad h_0= D_1=
c_0. \lab{hD} \ee The recurrence coefficients $u_n$ satisfy the
relation \be u_n = \frac{h_n}{h_{n-1}} = \frac{D_{n-1}
D_{n+1}}{D^2_n}. \lab{uh} \ee Thus we have \be h_n = c_0 u_1 u_2
\cdots u_n. \lab{hu} \ee
\par
Assume now that the polynomials $P_n(x;t)$ depend on a real
parameter $t$ through their recurrence coefficients $u_n(t)$,
$b_n(t)$. Then the restricted Toda chain equations (RTE) are
equivalent to the condition

\be \dot P_n(x;t) = - u_n P_{n-1}(x;t). \lab{dotP} \ee

It is possible to choose initial moment $c_0(t)$ (normalization)
such that  the RTE are equivalent to the very simple condition \be
\dot c_n = c_{n+1}, \lab{cnn} \ee i.e. \be c_n(t) = \frac{d ^n
c_0(t)}{dt^n}. \lab{cn0} \ee Hence, for the Toda chain case, the
Hankel determinants $D_n=D_n(t)$ have the form \be D_n(t) =
\det(c_{0}^{(i+k)}(t))_{i,k=0,\dots,n-1}, \quad D_0 =1, \quad D_1=
c_0, \lab{DeltaT} \ee where $c_{0}^{(j)}$ means the $j$-th
derivative of $c_0(t)$ with respect to $t$.

Under this condition, the RTE are equivalent also to the equations
\be \frac{d^2 \log D_n}{dt^2} = \frac{D_{n-1} D_{n+1}}{D_n^2},
\quad n=1,2,\dots. \lab{Sylv} \ee

Note also that for the Hankel determinants of the form \re{DeltaT}
we have the useful relation \be b_n = \frac{\dot D_{n+1}}{D_{n+1}}
- \frac{\dot D_n}{D_n} \lab{bD} \ee or, equivalently,  \be  b_n =
\dot h_n/ h_n. \lab{hb} \ee In particular, for $n=0$ we have from
\re{hb} \be b_0 = \frac{\dot c_0}{c_0}. \lab{b0c0} \ee The
relation \re{b0c0} allows us to restore $c_0(t)$ if the recurrence
coefficient $b_0=b_0(t)$ is known explicitly  from Toda chain
solutions \re{Toda}.

The Stieltjes function $F(z)$ is defined as a generating function
of the moments \cite{Chi} \be F(z) = \frac{c_0}{z }+
\frac{c_1}{z^2 }+ \cdots + \frac{c_n}{z^{n+1}}+ \cdots. \lab{Fz}
\ee If moments $c_n$ depend on $t$ according to the Toda Ansatz
\re{cnn}, we then have \be \dot F(z;t)= \frac{c_1 }{z} +
\frac{c_2}{z^2} + \cdots + \frac{c_n}{z^n} + \cdots = zF(z)-c_0.
\lab{dotF} \ee In fact, the relation \re{dotF} is equivalent to
the  restricted TC equations \re{cnn}.
\par
We consider also a so-called $E$-generating function of another
type: \be \Phi(p) = \sum_{k=0}^{\infty} c_k \frac{p^k}{k!}.
\lab{Phi} \ee The relationship between functions $F(z)$ and
$\Phi(p)$ is given by the (formal) Laplace transform: \be F(z)=
\sum_{k=0}^{\infty} c_k z^{-k-1} = \sum_{k=0}^{\infty} c_k
\int_{0}^{\infty} \frac{p^k e^{-pz}}{k!} dp = \int_{0}^{\infty}
e^{-pz} \Phi(p) dp. \lab{Lap} \ee

For the case of the RTE with condition \re{cnn} we see that
generating function $\Phi(p)$ is given automatically by the formal
Taylor expansion \be \Phi(p;t)= \sum_{k=0}^{\infty} c_k(t)
\frac{p^k}{k! }= \sum_{k=0}^{\infty} c_0^{(k)}(t) \frac{p^k}{k! }=
c_0(t+p) \lab{Phi1} \ee of $c_0(t+p)$. Thus the $E$-generating
function is given just by the shifted $c_0(t+p)$ zero-moment
function. The Stieltjes function is given then as the Laplace
transform \be F(z;t) =\int_{0}^{\infty} e^{-pz} c_0(t+p) dp.
\lab{Fz1} \ee It should be noted, however, that formula \re{Fz1}
has rather formal meaning. In practice, there are situations when
direct application of this formula may be problematic, if, e.g.
integral in rhs \re{Fz1} diverges. As a simple example consider
the case when $c_0(t)$ has the expression \be c_0(t) =
\sum_{k=-\infty}^{\infty} \mu_k \exp(\nu_k t) \lab{c_0_exp} \ee
with some complex constants $\mu_k, \nu_k$. Then
$$
c_n(t) = \frac{d^n c_0(t)}{dt^n} = \sum_{k=-\infty}^{\infty} \mu_k \nu_k ^n \exp(\nu_k t)
$$
and the Stieltjes function $F(z;t)$ has the expression \be F(z;t)
= \sum_{n=0}^{\infty} c_n(t) z^{-n-1} = \sum_{k=-\infty}^{\infty}
\frac{\mu_k \exp(\nu_k t)}{z-\nu_k} \lab{F(z)_exp} \ee i.e. the
orthogonality measure in this case is located at the points
$\nu_k$ of the complex plane with the corresponding concentrated
masses $M_k(t) = \mu_k \exp(\nu_k t)$. Formula \re{F(z)_exp} is a
special case of the well known result for the restricted Toda
chain: if $d \rho(x)$ is an orthogonality measure for the
orthogonal polynomials $P_n(z;0)$ (i.e. for initial value of time
$t=0$), then for arbitrary $t$ the measure will be \cite{Apt},
\cite{Peh1} \be d \rho(x;t) = const \: \exp(xt) \: d \rho(x)
\lab{rho_t} \ee where the constant in rhs of \re{rho_t} is not
essential and is needed only to provide the normalization
condition for the measure.  Formula \re{F(z)_exp} corresponds to
the special case of the measure \be d\rho(x) =
\sum_{k=-\infty}^{\infty} \mu_k \: \delta(x-\nu_k)dx
\lab{rho_discr} \ee On the other hand, the generating function
$\Phi(p;t)$ has the expression \be \Phi(p;t) = \sum_{n=0}^{\infty}
\frac{c_n(t) p^n}{n!} = \sum_{k=-\infty}^{\infty} \mu_k \exp(\nu_k
(p+t)) \lab{Phi_exp} \ee If one applies formula \re{Lap} for the
Laplace transform to the function $\Phi(p;t)$ given by
\re{Phi_exp} we get \be F(z;t) = \sum_{k=-\infty}^{\infty}e^{\nu_k
t} \mu_k \int_0^{\infty} e^{-p(z-\nu_k)} dp \lab{Lap1} \ee Doing
"naively" we can put \be \int_0^{\infty} e^{-p(z-\nu_k)} dp =
(z-\nu_k)^{-1} \lab{Lap2} \ee in \re{Lap1} and obtain desired
formula \re{F(z)_exp} for the Stieltjes function $F(z;t)$.
However, formula \re{Lap2} is correct only if $Re(z)>Re(\nu_k)$.
For $Re(z)<Re(\nu_k)$ we may still use the same formula \re{Lap1}
because it corresponds to the true formal series \re{F(z)_exp}.

Sometimes the following trick will be useful.   Consider the
modified function $\tilde c_0(t) = c_0(it)$ (i.e. just pass to the
"imaginary time"). Construct the Hankel determinants $\tilde
D_n(t)$ form the new moments $\tilde c_n(t) = i^n c_n(it)$.
Clearly, $\tilde D_n(t)= i^{n(n-1)} D_n(it)$. Define corresponding
recurrence coefficients $\tilde b_n(t), \tilde u_n(t)$. We have
\be \tilde b_n(t) = i b_n(it), \quad \tilde u_n(t) = - u_n(it)
\lab{mod_bu} \ee Obviously they will satisfy the restricted Toda
chain equations \re{Toda}. Construct corresponding orthogonal
polynomials \be \tilde P_n(z;t) = i^{n}P_n(z/i;it) \lab{mod_Pn}
\ee These "modified" orthogonal polynomials satisfy the recurrence
relation \be \tilde P_{n+1}(z;t) + \tilde b_n(t) \tilde P_n(z;t) +
\tilde u_n \tilde P_{n-1}(z;t) = z \tilde P_n(z;t) \lab{rec_mod_P}
\ee with the initial conditions
$$
\tilde P_0(z;t) =1, \quad \tilde P_1(z;t) = z- \tilde b_0(t)
$$
The Stieltjes function $\tilde F(z;t)$ for these polynomials is
defined as
$$
\tilde F(z;t) = \sum_{n=0}^{\infty} \tilde c_n(t) z^{-n-1} =
\sum_{n=0}^{\infty} i^n c_n(it) z^{-n-1} =-i F(z/i,it)
$$
Assume now that the function $\tilde c_0(t)$ is periodic with some
real period $\tilde c_0(t+T)=\tilde c_0(t)$. We assume also that
the function $\tilde c_0(t)$ is regular in some maximal strip
$-\alpha < y < \beta$, where $y = \Im(t)$ and $\alpha, \beta$ are
some positive parameters. It is well known (see, e.g.
\cite{Akhiezer}) that inside this strip the function $\tilde
c_0(t)$ has no singularities and it has at least one singularity
point on both lines $y=-\alpha$ and $y = \beta$.

Inside this strip we have the Fourier series expansion \be \tilde
c_0(t) = \sum_{k=-\infty}^{\infty} \mu_k \exp(2\pi i k t/T)
\lab{Four_c0} \ee For the Fourier coefficients there is the useful
asymptotic estimation \cite{Akhiezer} \be {  |\mu_{-k}| \le Q\:
\exp(\frac{-2\pi k \beta}{T}), \atop |\mu_{k}| \le Q \:
\exp(\frac{-2\pi k \alpha}{T})} \lab{strip} \ee for all
$k=0,1,2,\dots$ with a constant $Q$ not depending on $k$.

In this case the measure for the modified polynomials $\tilde
P_n(z;t)$ is purely discrete and is located on the uniform grid
\be z_k = 2 \pi i k /T \lab{im_grid}, \quad k=0, \pm 1, \pm
2,\dots \ee on the imaginary axis. Hence, for real and periodic
functions $\tilde c_0(t)$ we obtain orthogonal polynomials $\tilde
P_n(z;t)$ with the measure located on the imaginary axis.

Returning to the initial orthogonal polynomials $P_n(z;t)$
we have the series
\be
c_0(t) =
\sum_{k=-\infty}^{\infty} \mu_k \exp(2\pi k t/T) \lab{exp_c0} \ee
which for real values of the argument $t$ converges on the interval $-\beta< t < \alpha$.

We thus obtain \be F(z;t) = \sum_{k=-\infty}^{\infty} \frac{\mu_k
}{z- 2 \pi k /T} \: \exp(2 \pi kt/T) \lab{FF_im} \ee

We see that in this case the Stieltjes function $F(z;t)$
corresponds to a purely discrete measure located on the {\it real}
axis at points \be z_k =  2\pi k/T, \; k=0, \pm 1 , \pm 2, \dots
\lab{real_grid} \ee with the corresponding discrete masses \be
M_k(t) = \mu_k \: \exp(2 \pi k t/T) \lab{mod_M_k} \ee The measure
is well defined for all real values of the parameter $t$ belonging to the admissible interval $-\beta< t < \alpha$.
Indeed we have the expression for the moments \be c_j(t)=
\sum_{k=-\infty}^{\infty} M_k(t)z_k^j = \sum_{k=-\infty}^{\infty}
\mu_k \: \exp(2 \pi k t/T) z_k^j. \lab{restr_M} \ee From estimations \re{strip} it follows that for all $j=0,1,2,\dots$
the corresponding sum in \re{restr_M} converges in the interval $-\beta< t < \alpha$ and hence all the moments $c_j(t)$ are well defined
for $t$ belonging to this admissible interval.

The remaining problem is to determine when the moment problem is
positive definite, i.e. we need to find when the condition $D_n(t)
>0$ holds for all $n=0,1,2,\dots$ and for all $t$ belonging to the
admissible interval $-\beta< t < \alpha$. It is well known from
general theory of orthogonal polynomials \cite{Ismail}, \cite{Chi}
that in our case this condition is equivalent to the positivity of
all concentrated masses $M_k(t) >0, \; k=0, \pm 1, \pm 2, \dots$
if $t$ belongs to the admissible interval. Obviously, this
condition is equivalent to positivity of all Fourier coefficients
$\mu_k >0, \; k=0, \pm 1, \pm 2, \dots$ for the periodic function
$\t c_0(t)$ when $t$ belongs to the maximal strip $-\alpha <
\Im(t) < \beta$.

In turn, this condition means that the function $\t c_0(t)$
belongs to the class of so-called positive-definite functions.

Recall \cite{Akh_mom} that the continuous function $f(x)$ is
called the positive-definite if  for any n, any real parameters
$x_i, \; i=1,2,\dots, n$ and any complex variables $\xi_i, \;
i=1,2,\dots n$ the property \be \sum_{i,k=1}^n f(x_i-x_k)\xi_i
\bar \xi_k \ge 0 \lab{Hdf} \ee holds. (In the probability theory
the positive-definite functions are called the characteristic
functions \cite{Lukacs}). This property is equivalent to a
possibility to present the function $f(x)$ in the form \be f(x) =
\int_{-\infty}^{\infty} e^{ixt} d \sigma(t), \quad -\infty < x<
\infty, \lab{Boch_f} \ee where $\sigma(t)$ is a nondecreasing
function of a bounded variation. For the proof of this important
proposition (the Bochner theorem) see, e.g. \cite{Akh_mom},
\cite{Lukacs}.

From this property it is elementary derived that the function is
bounded $|f(x)|\le f(0)$ and $f(0)>0$. Moreover \be f(-x) = {\bar
f(x)} \lab{T_f} \ee

A special case is positive definite {\it periodic} functions.
Assume that $f(x)$ is periodic with some real period $T$:
$f(x+T)=f(x)$. We have the Fourier expansion \be f(x) =
\sum_{n=-\infty}^{\infty} A_n \: e^{2\pi i n x/T}, \lab{Four_f}
\ee The periodic function $f(x)$ is positive definite if and only
if all the Fourier coefficients are nonnegative $A_n \ge 0$.
Moreover, we assume that there are infinity many positive
coefficients $A_n >0$ (otherwise the problem is trivial).

In our case we have that the periodic function $\t c_0(t)$ should
have only nonnegative Fourier coefficients $\mu_k$.

We thus have
\begin{pr}
Assume that $\t c_0(t) = c_0(it)$ is a periodic function with a
real period $T$ and a non-empty maximal strip of regularity
$-\alpha < \Im(t) < \beta$ with $\alpha, \beta$ some positive
constants. Then the moment problem corresponding to the function
$c_0(t)$ will be positive definite (i.e. $D_n(t)>0, \:
n=0,1,2,\dots$ when $-\beta < t <\alpha$) if and only if the
function $\t c_0(t)$ is positive definite.
\end{pr}
Note that in \cite{Zhe_cndn} we considered explicit examples of
polynomials orthogonal on the unit circle connected with positive
definite functions. In that case positive definite Toeplitz
determinants (they correspond to a positive nondecreasing measure
on the unit circle) appear instead of the Hankel ones.

We see that the trick with passing to the purely imaginary time $t
\to it$ is especially useful if the Fourier series for the
modified function $\tilde c_0(t)$ is known. Then we can restore
information on the initial measure using formula \re{FF_im}.

Let $F(z;t)$ and $G(z;t)$ be the Stieltjes functions corresponding
to the moments $c_0(t)$ and $g_0(t)= e^{\gamma t} c_0(t)$, where
$\gamma$ is a complex constant.  Then it follows from \re{Lap}
that \be G(z;t) = e^{\gamma t} F(z-\gamma;t). \lab{g_F} \ee This
means that if $P_n(x;t)$ are monic orthogonal polynomials
corresponding to the moments $c_n(t)=\frac{d^n c_0(t)}{dt^n}$ then
$P_n(x-\gamma;t)$ are monic orthogonal polynomials corresponding
to the moments $g_n(t)=\frac{d^n g_0(t)}{dt^n}$.

This simple observation may be useful e.g. in the case when $\t
c_0(t) = c_0(it)$ is {\it quasi-periodic}, i.e. it satisfies the
condition \be \t c_0(t+T) = e^{\nu} \t c_0(t) \lab{quasi_c0} \ee
with a real period $T$ and some complex constant $\nu$. Then we
can introduce the new moment $\t g_0(t) = \exp(-\nu t/T) \t
c_0(t)$ which is purely periodic: $\t g_0(t+T) = \t g_0(t)$.
Assume that the function $\t g_0(t)$ is positive definite and has
a nonempty regularity strip in the complex domain. Starting with
the moment $g_0(t)= \t g_0(-it)$ one can construct monic
orthogonal polynomials $Q_n(x;t)$ having purely discrete positive
measure on the real line. Then the orthogonal polynomials
$P_n(x;t)$ corresponding to the moment $c_0(t)$ have the
expression $P_n(x;t)= Q_n(x + i \nu/T;t)$. The measure for the
polynomials $P_n(x;t)$ will be located on the horizontal line
$\Im(z) = \nu_1/T$, where $\nu_1 = \Re(\nu)$. Thus the
quasi-periodic functions $\t c_0(t)$ lead to a simple shift of the
argument $x \to x+ i \nu/T$ of orthogonal polynomials $P_n(x;t)$
corresponding to a pure periodic function $\t c_0(t)$.

\section{Elliptic functions.}
\setcounter{equation}{0} In this section we recall basic
properties of the Weierstrass elliptic functions which will be
needed in further analysis \cite{Akhiezer}, \cite{WW}. The
Weierstrass function $\wp(z;g_2,g_3)$ depends on the argument $z$
and two parameters $g_2,g_3$ (the so-called invariants). It
satisfies the differential equation
$$
\wp'^2(z) = 4 \wp^3(z) -g_2 \wp(z) -g_3 = 4(\wp(z) - e_1)(\wp(z) - e_3)(\wp(z) - e_3)
$$
where the parameters $e_i$ satisfy the restriction $e_1+e_2+e_3=0$. The function $\wp(z)$ is double-periodic:
$$
\wp(z+2\omega) =\wp(z+2\omega') = \wp(z),
$$
where the periods $2 \omega, 2\omega'$ are assumed to satisfy the
condition $Im(\omega'/\omega) >0$. If $g_2, g_3$ are known then
the periods $2 \omega, 2\omega'$ can be calculated by a standard
procedure in terms of the elliptic integrals of the first kind
\cite{Akhiezer}.

It is convenient to introduce the notation \cite{Akhiezer}
$$
\omega_1 = \omega, \; \omega_3 = \omega', \; \omega_2 = -\omega - \omega'
$$
There is a relation between $e_i$ and $\omega_i$:
\be
\wp(\omega_k) = e_k, \quad k=1,2,3 \lab{wek} \ee
The Weierstrass zeta function $\zeta(z)$ is an odd function $\zeta(-z) = - \zeta(z)$ defined as
$$
\zeta'(z) = - \wp(z)
$$
The function $\zeta(z)$ has simple poles at the points $2m\omega + 2m'\omega'$, where $m, m'$ are arbitrary integers.
In contrast to $\wp(z)$, the zeta function $\zeta(z)$ is quasiperiodic:
\be
\zeta(z+ 2 \omega_k) = \zeta(z) + 2 \eta_k, \quad k=1,2,3
\lab{per_z} \ee
where
$$
\eta_k = \zeta(\omega_k)
$$
There are useful relations for $\eta_k$:
\be
\eta_1 + \eta_2 + \eta_3 =0 \lab{eee} \ee
and
$$
\eta_2 \omega_1 - \eta_1 \omega_2 = i \pi/2
$$
(there are two similar relations which are obtained from the last
relation by a cyclic permutation of 1,2,3).

The Weierstrass sigma function $\sigma(z)$ is an odd function
$\sigma(-z)=-\sigma(z)$ defined as \be
\frac{\sigma'(z)}{\sigma(z)} = \zeta(z) \lab{sigma_zeta} \ee It
has simple zeroes at the points $2m\omega + 2m'\omega'$. The sigma
function has quasi-periodic property \be \sigma(z+2\omega_k) =
-\exp(2\eta_k (z+\omega_k)) \: \sigma(z) \lab{per_s} \ee There is
a simple formula connecting $\wp(z)$ and $\sigma(z)$
\cite{Akhiezer}: \be \frac{\sigma(u+v) \sigma(u-v)}{\sigma^2(u)
\sigma^2(v)} = \wp(v) - \wp(u) \lab{s_wp} \ee Apart from the
function $\sigma(z)$ one can define functions $\sigma_k(z),
k=1,2,3$ by the formulas \cite{Akhiezer} \be \sigma_{\alpha}(z) =
\frac{\sigma(z+\omega_{\alpha})}{\sigma(\omega_{\alpha})}\:
\exp(-z\eta_{\alpha}), \quad \alpha=1,2,3 \lab{sig_k} \ee The
functions $\sigma_{\alpha}(z)$ are convenient when passing from
the Weierstrass to the Jacobi elliptic functions. Indeed, we have
\cite{WW} \be \sn(u;k) = (e_1-e_3)^{1/2}
\frac{\sigma(z)}{\sigma_3(z)}, \; \cn(u;k) =
\frac{\sigma_1(z)}{\sigma_3(z)}, \; \dn(u;k) =
\frac{\sigma_2(z)}{\sigma_3(z)} \lab{Jac_sig} \ee where \be
u=(e_1-e_3)^{1/2}z, \quad k^2= \frac{e_2-e_3}{e_1-e_3} \lab{uzk}
\ee The parameter $k$ is called the elliptic modulus. The
parameter
$$
k' = (1-k^2)^{1/2} = \sqrt{\frac{e_1-e_2}{e_1-e_3}}
$$
is called the complementary modulus \cite{WW}. The values
$$
K = \sqrt{e_1-e_3} \: \omega_1, \quad K' = i \sqrt{e_1-e_3} \:
\omega_3
$$
are complete elliptic integrals of the first kind \cite{WW}.

We need also expressions of the Weierstrass functions for the
value $z=\omega_1/2$: \ba &&\wp(\omega_1/2) = \wp(3 \omega_1/2) =
e_1 + (e_1-e_3)k', \; \wp(\omega_3 + \omega_1/2) = \frac{e_3k' +
e_2}{1+k'}, \nonumber \\ && \wp'(\omega_1/2) = - 2 k'(1+k')
(e_1-e_3)^{3/2}, \; \wp''(\omega_1/2) = 4(e_1-e_3)
\left(2(e_1-e_2)  + 3 e_1 k' \right) \lab{wp_half} \ea \ba &&
2\zeta(\omega_1/2) = \eta_1 - \frac{1}{2} \:
\frac{\wp''(\omega_1/2)}{\wp'(\omega_1/2)} = \nonumber \\&& \eta_1
+ \sqrt{e_1-e_3} (k'+1), \nonumber
\\&& \zeta(\omega_3 + \omega_1/2) = \eta_3 + \eta_1/2 + \sqrt{e_1-e_3} (1-k')/2. \lab{zeta_half} \ea
Note that the choice of an appropriate sign in front of square
roots in formulas \re{wp_half} and \re{zeta_half} is not trivial
problem and depends on location of the parameters $e_1,e_2,e_3$ in
the complex domain. However in our further analysis we will use
the "canonical" choice of these parameters: they are real and
ordered as $e_3 <e_2<e_1$. Then all the square roots are assumed
in the arithmetic meaning.

Apart from the Weierstrass zeta function $\zeta(z;g_2,g_3)$
sometimes the Jacobi Zeta function $Z(z;k)$ is more convenient.
The relation between these function is \cite{WW} \be \zeta(z) =
\frac{z\eta_1}{\omega_1} + \sqrt{e_1-e_3}\left\{ Z(u,k) +
\frac{\cn(u,k) \dn(u,k)}{\sn(u,k)}\right\}, \lab{zZ} \ee where the
same relations \re{uzk} are assumed.

The Jacobi Zeta function is purely periodic with respect to the
period $2 K$:
$$
Z(u+2K) = Z(u)
$$
and quasi-periodic with respect to the period $2iK'$:
$$
Z(u+2iK') = Z(u) - \frac{i\pi}{K}
$$
It possesses a remarkable "addition theorem" \cite{WW} \be Z(u+v)
= Z(u) + Z(v) - k^2 \sn(u)\sn(v) \sn (u+v) \lab{add_Z} \ee Note
the useful relations \be Z(u+K)= Z(u) - k^2 \: \frac{\sn(u)
\cn(u)}{\dn(u)} \lab{Z_K} \ee \be Z(u+iK')= Z(u) +  \frac{\dn(u)
\cn(u)}{\sn(u)}  - \frac{i \pi}{2K}  \lab{Z_K1} \ee \be
Z(u+K+iK')= Z(u) - \frac{\sn(u) \dn(u)}{\cn(u)}  - \frac{i
\pi}{2K}.  \lab{Z_KK1} \ee The Jacobi Zeta function can be
expressed in terms of the incomplete elliptic integral $E(u)$ of
the second kind \cite{WW} \be Z(u) = E(u) - u \: \frac{E}{K}
\lab{ZE} \ee where
$$
E(u) = \int_{0}^u
 \dn^2(t) dt $$ and $E= E(K)$.
The following formulas allow to express the functions
$\zeta(z+\omega_i), \: i=1,2,3$ in terms of the Jacobi Zeta
function: \be \zeta(z+\omega_1) =\eta_1+ \frac{z\eta_1}{\omega_1}
+ \sqrt{e_1-e_3}\left\{ Z(u,k) - \frac{\sn(u,k)
\dn(u,k)}{\cn(u,k)}\right\}, \lab{zZ1} \ee

\be \zeta(z+\omega_3) =\eta_3 + \frac{z\eta_1}{\omega_1} +
\sqrt{e_1-e_3} Z(u,k)
  , \lab{zZ3} \ee

\be \zeta(z+\omega_1+\omega_3) =-\eta_2+ \frac{z\eta_1}{\omega_1}
+ \sqrt{e_1-e_3}\left\{ Z(u,k) - k^2\: \frac{\sn(u,k)
\cn(u,k)}{\dn(u,k)}\right\}, \lab{zZ31} \ee

\section{Toda chain solution and the corresponding orthogonal
polynomials} \setcounter{equation}{0} We present here an explicit
solution of the restricted Toda chain
\begin{lem}
Put \be u_n(t) = w^2n^2 \left(\wp(w(t+\beta)) - \wp(nw(t+\beta)
+q) \right) \lab{u_FS} \ee and \be b_n(t) = \mu_1 + w(n+1)
\zeta(w(n+1)(t+\beta) + q) - wn \zeta(wn(t+\beta) + q) -(2n+1) w
\zeta(w(t+\beta)), \lab{b_FS} \ee where $w,\beta,q, \mu_1$ are
arbitrary complex parameters. We will also assume that $\omega_1,
\omega_3$ are arbitrary independent periods corresponding to the
arbitrary parameters $e_1,e_2,e_3$ with the only condition
$e_1+e_2+e_3=0$.

Then $u_n(t), b_n(t)$ satisfy the restricted Toda chain equations
\re{Toda}
\end{lem}

{\it Proof}. In order to verify the first equation in \re{Toda} we
present $u_n(t)$ in an equivalent  form \be u_n(t) = w^2 n^2 \:
\frac{\sigma((n+1) w(t+\beta) +q) \sigma((n-1) w(t+\beta) +q)
}{\sigma^2(nw(t+\beta) +q) \sigma^2(w(t+\beta)) } \lab{U_sigma}
\ee using formula \re{s_wp}. Then by \re{sigma_zeta} the
expression $\dot u_n/u_n$ can be presented as  a sum of the
Weierstrass zeta functions and we arrive at the first equation of
\re{Toda}. The second equation \re{Toda} is satisfied by the
formula $\wp(z) = -\zeta'(z)$.

It is directly verified that the recurrence coefficients $u_n(t),
b_n(t)$ are double-periodic with the periods $2\omega/w, 2
\omega'/w$:
$$
u_n(t+2\omega/w)=u_n(t+2\omega'/w) = u_n(t), \quad b_n(t+2\omega/w)=b_n(t+2\omega'/w) = b_n(t)
$$
(periodicity property for $u_n(t)$ is obvious and periodicity for
the coefficients $b_n(t)$ follows from \re{per_z} and \re{eee}).
Thus both $u_n(t)$ and $b_n(t)$ are elliptic functions in the
argument $t$.

Using formulas \re{b0c0} and \re{sigma_zeta} we can restore the
function $c_0(t)$. It is easy to verify that \be c_0(t) =
\frac{\sigma(w(t+\beta) + q)}{\sigma(q) \: \sigma(w(t+\beta))} \:
\exp(\mu_1 t + \mu_0), \lab{c0_FS} \ee where $\mu_0$ is an
arbitrary constant.

We thus obtained some new family of orthogonal polynomials
$P_n(z;t)$ which can be defined through given recurrence
coefficients $u_n(t), b_n(t)$. The Stieltjes function $F(z,t)$
(and hence, in principle) the orthogonality measure for these
polynomials can also be found explicitly from formula \re{Fz1}
because the function $c_0(t)$ is given explicitly by \re{c0_FS}.

The obtained orthogonal polynomials $P_n(z;t)$ contain several
free parameters (say $w,q,\beta, \mu_1, t$ and elliptic parameters
$g_2, g_3$). We would like to investigate some simple special
choice of these parameters when our polynomials are a
generalization of already known families. Note that the parameter
$\beta$ is inessential: it describes a shift of the argument $t
\to t+\beta$. Nevertheless, we will keep this parameter for
convenience, assuming that the argument $t$ takes real or pure
imaginary  values. We will assume also that $q \ne 0$. Indeed, the
case $q=0$ corresponds to some degeneration: $b_0(t) =const, \;
u_1(t) =0$ and $c_0(t)$ becomes a pure exponential function $c_0 =
\exp(\mu t)$ in this limit.

In what follows we put $q=\omega_j, \: \beta = \omega_k/w$, where
$k,j$ are arbitrary noncoinciding integers from the set $1,2,3$.
We denote also $\omega_l =-\omega_j-\omega_k$. Then using
(quasi)periodicity properties of the Weierstrass functions
$\wp(z), \zeta(z)$ we find the expression for $u_n(t)$: \ba
&&u_{2n}(t) = 4
w^2 n^2 (\wp(wt + \omega_k) - \wp(2wnt + \omega_j)) \nonumber \\
&& u_{2n+1}(t) =  w^2 (2n+1)^2 (\wp(wt + \omega_k) - \wp(w(2n+1)t
+ \omega_l)) \lab{u_om} \ea and for the coefficients $b_n(t)$: \ba
&& b_{2n}(t) = \mu_1 + w\left\{(2n+1) \zeta((2n+1)wt - \omega_l) -
2n \zeta(2nwt+\omega_j) - (4n+1)\zeta(wt+\omega_k) + 2n \eta_k
\right\}, \nonumber \\&& b_{2n+1}(t) = \mu_1 + w\left\{(2(n+1)
\zeta((2(n+1)wt + \omega_j) - (2n+1) \zeta((2n+1)wt-\omega_l) - \right . \nonumber \\
&& \left . (4n+3)\zeta(wt+\omega_k) + (6n+4) \eta_k     \right\} \lab{b_om} \ea

We will assume also that $\mu_0=0$ and \be \mu_1 = -w \eta_j,
\lab{mu_1_eta} \ee Indeed, the parameter $\mu_0$ is inessential
and can be chosen arbitrary whereas the parameter $\mu_1$ leads
only to a trivial shift of the recurrence coefficient $b_n$, hence
we can put $\mu_1$ to a prescribed value without loss of
generality.

In order to find the orthogonality measure for the obtained
polynomials $P_n(z;t)$ we will assume that the parameter $w$ is
real. Among all 6 possible choices $q=\omega_j, \: \beta =
\omega_k/w$ of the parameters $q,\beta$ we consider only the two
cases:

(i) if $\beta= \omega_1/w,  \; q = \omega_2$ then form \re{sig_k},
\re{mu_1_eta}  and \re{Jac_sig} we find that \be c_0(t) =
C_1/\cn(w\sqrt{e_1-e_3} t;k) \lab{c_0_cn} \ee where
$$
C_1 = -\frac{\sigma(\omega_3)}{\sigma(\omega_1)\sigma(\omega_2)}\:
e^{-\omega_1 \eta_2}
$$

(ii) if $\beta= \omega_1/w,  \; q = \omega_3$ then quite
analogously we find \be c_0(t) = C_2 \: \frac{\dn(w\sqrt{e_1-e_3}
t;k)}{\cn(w\sqrt{e_1-e_3} t;k)} = C_2 \: \dc(w\sqrt{e_1-e_3}
t;k)\lab{c_0_dn} \ee where
$$
C_2 = -\frac{\sigma(\omega_2)}{\sigma(\omega_1)\sigma(\omega_3)}\:
e^{-\omega_1 \eta_3}
$$

 Note that the constant factors $C_1, C_2$
are in fact inessential (the orthogonal polynomials $P_n(z;t)$ as
well as recurrence coefficients $b_n(t), u_n(t)$ do not depend on
these constants) and we can put $C_2=C_1=1$ . We thus have that
$c_0(t) = 1/\cn(wt\sqrt{e_1-e_3};k)$ for the case (i) and $c_0(t)
= \dc(w\sqrt{e_1-e_3} t;k))$ for the case (ii).

We will restrict ourselves with the case when all $e_i$ are real
and distinct, say $e_1>e_2>e_3$ (this is the usual convention
\cite{Akhiezer}). Then it is well known \cite{Akhiezer} that  for
the real values of $z$ the functions $\wp(z), \zeta(z), \sigma(z)$
take real values. For purely imaginary values of $z$ the function
$\wp(z)$ takes real values whereas functions $\zeta(z), \sigma(z)$
take purely imaginary values. The period $2 \omega_1$ is real and
the period $2 \omega_3$ is purely imaginary. This means that the
fundamental parallelogram of the elliptic functions in this case
is a rectangle \cite{Akhiezer}. We have also that both $k$ and
$k'$ are real parameters taking values in the "canonical" interval
$0< k,k' < 1$. Hence all values of Jacobi elliptic functions
$\sn(x),\cn(x), \dn(x)$ are real for real $x$. The functions
$\cn(x), \dn(x)$ take also real values for purely imaginary values
of $x$.

Now we are ready to  calculate the orthogonality measure for the
cases (i) and (ii).

\section{The orthogonality measure and recurrence coefficients for the case (i)}
\setcounter{equation}{0} According to considerations of the first
section, introduce the new function $\tilde c_0(t) = c_0(it)$.

For the first case (i) we have \be \tilde c_0(t) = 1/\cn(i
w\sqrt{e_1-e_3}t;k) = \cn(w\sqrt{e_1-e_3}t;k') \lab{real_cn} \ee

The function $\cn(w\sqrt{e_1-e_3}t;k')$ is real on the real axis
and periodic with the period $T=\frac{4K'}{w\sqrt{e_1-e_3}}$. The
Fourier series for this function is well known \cite{WW} \be \cn(w
\sqrt{e_1-e_3}t; k') = \frac{\pi }{k' K'}
\sum_{n=-\infty}^{\infty} \frac{1}{v^{n-1/2} + v^{1/2-n} }
\exp(\pi i (n-1/2)w \sqrt{e_1-e_3} t /K'), \lab{cn_F} \ee where
$$
v= \exp(-\pi K/K')
$$
Hence the polynomials $P_n(z;t)$ have a purely discrete
orthogonality measure located at the points \be x_n = \frac{2
\pi}{T}\: (2n-1) = \frac{\pi w \sqrt{e_1-e_3}}{2 K'} (2n-1), \quad
n=0, \pm 1, \pm 2 , \dots \lab{loc_cn} \ee with the corresponding
concentrated masses \be M_n(t) = \frac{\pi}{k'K'} \:
\frac{\exp(\pi w t(n-1/2) \sqrt{e_1-e_3} /K')}{v^{n-1/2} +
v^{1/2-n}} \lab{M_cn} \ee Thus orthogonality relation for the
polynomials $P_n(z;t)$ looks as follows \be
\sum_{s=-\infty}^{\infty} {M_s(t) P_n(x_s;t)P_m(x_s;t)} = h_n(t)
\: \delta_{nm} \lab{ort_cn} \ee It is interesting to determine
conditions under which the measure is well defined, i.e. that all
the moments are finite \be \tilde c_j(t) =
\sum_{s=-\infty}^{\infty} M_s(t) x_s^j < \infty, \quad
j=0,1,2,\dots \lab{M_cond} \ee From the explicit expression
\re{M_cn} it is easily seen that condition \re{M_cond} will hold
provided \be -\frac{K}{w \sqrt{e_1-e_3}}< t < \frac{K}{w
\sqrt{e_1-e_3}} \lab{t_ad_cn} \ee  When the parameter $t$ belongs
to this interval all the moments are well defined. If $t \to \pm
\frac{K}{w \sqrt{e_1-e_3}}$ then $c_0(t) \to \infty$ as is easily
seen from \re{real_cn}. Hence when $t$ approaches the endpoints of
the interval \re{t_ad_cn}, the moments $c_n(t)$ tend to infinity
and the measure becomes not well defined.

Note that the admissible interval \re{t_ad_cn} corresponds to the strip of regularity
for of the function $\tilde c_0(t)$. In turn, this strip of regularity is obtained from the strip
of convergence of the Fourier series for the elliptic Jacobi function $\cn(z;k)$ \cite{Akhiezer}.

From \re{M_cn} it is
clear that for all values of $t$ from the admissible interval
\re{t_ad_cn} the concentrated masses are positive $M_s(t)>0, \;
s=0,\pm 1, \pm 2, \dots$. This means that we indeed deal with a
positively defined purely discrete measure on the whole real axis.

Moreover it is easy verified that \be \sum_{s=-\infty}^{\infty}
M_s(t) x_s^j = c_j(t) = \frac{d^j}{dt^j} \: c_0(t) \lab{true_moms}
\ee where $c_0(t)=1/\cn(wt\sqrt{e_1-e_3};k)$. Indeed, formula
\re{true_moms} follows directly from the Fourier series \re{cn_F}
by $j$-fold differentiation with respect to $t$. Formula
\re{true_moms} shows that the obtained measure is "true", i.e. it
gives the prescribed moments $c_j(t)$ for all $j=0,1, \dots$.

Consider the recurrence coefficients $b_n(t), u_n(t)$ for the
orthogonal polynomials $P_n(z;t)$ corresponding to the function
$\tilde c_0(t)$ defined by \re{real_cn}.

It is clear that both $b_n(t)$ and $u_n(t)$ are real for all $t$
from the admissible interval \re{t_ad_cn}. Indeed, the Hankel
determinants $D_n(t)$ are real because these are constructed from
the matrix with real entries
$$
a_{ij} =\frac{d^{i+j}  c_0(t)}{dt^{i+j}}, \; i,j=0,1,\dots, n-1
$$
Hence the normalization coefficients $h_n(t)$ are real as well.
The same is true for coefficients $b_n(t), u_n(t)$ obtained from
$h_n(t)$ by formulas
$$
b_n(t) = \frac{\dot h_n}{h_n}, \quad u_n = h_n/h_{n-1}
$$
Explicitly we have for the recurrence coefficients $b_n(t)$ (see
\re{b_om}) \ba && b_{2n}(t) =  w\left\{(2n+1) \zeta((2n+1)wt -
\omega_3) - 2n \zeta(2nwt+\omega_2) - (4n+1)\zeta(wt+\omega_1) +
2n \eta_1 - \eta_2 \right\}, \nonumber \\&& b_{2n+1}(t) =
w\left\{(2(n+1)
\zeta((2(n+1)wt + \omega_2) - (2n+1) \zeta((2n+1)wt-\omega_3) - \right . \nonumber \\
&& \left . (4n+3)\zeta(wt+\omega_1) + (6n+4) \eta_1   - \eta_2
\right\} \lab{b_cn} \ea

or, in terms of the elliptic Jacobi functions

\ba &&b_{2n}(t) = w \sqrt{e_1-e_3} \: \left\{ (2n+1)Z((2n+1)u) -
2n Z(2nu)  - (4n+1) Z(u) - \nonumber \right .\\ &&\left . 2n k^2
\: \frac{\cn(2nu) \sn(2nu)}{\dn(2nu)} + (4n+1) \: \frac{\sn(u)
\dn(u)}{\cn(u)} \right\} = \lab{b_even_Z} \\
&& w \sqrt{e_1-e_3} \: \left\{ (2n+1)Z((2n+1)u) - 2n Z(2nu+K)  -
(4n+1) \left(Z(u+K+iK') + \frac{i \pi}{2K} \right)\right\}
\nonumber \ea

\ba &&b_{2n+1}(t) = w \sqrt{e_1-e_3} \: \left\{ 2(n+1)Z(2(n+1)u) -
(2n+1) Z((2n+1)u)  - (4n+3) Z(u) - \nonumber \right .\\ &&\left .
2(n+1) k^2 \: \frac{\cn(2(n+1)u) \sn(2(n+1)u)}{\dn(2(n+1)u)} +
(4n+3) \: \frac{\sn(u) \dn(u)}{\cn(u)} \right\} = \lab{b_odd_Z}\\
&&w \sqrt{e_1-e_3} \: \left\{ 2(n+1)Z(2(n+1)u+K) - (2n+1)
Z((2n+1)u) - (4n+3) \left( Z(u) + \frac{i \pi}{2K} \right)
\right\} \nonumber \ea

where
$$
u= w \: t \: \sqrt{e_1-e_3}
$$

For the recurrence coefficients $u_n(t)$ we have expressions in
terms of the elliptic Jacobi functions \ba &&u_{2n}(t)= 4 n^2 w^2
(e_1-e_2)\: \left(\frac{1}{\cn^2(u)} + k^2 \:
\frac{\sn^2(2nu)}{\dn^2(2nu)} \right) \nonumber \\&& u_{2n+1}(t)=
(2n+1)^2 w^2 (e_1-e_3)\: \left(k'^2\: \frac{\sn^2(u)}{\cn^2(u)} +
\dn^2((2n+1)u) \right) \lab{u_cn} \ea

It is seen from \re{u_cn} that for all $t$ from the admissible
interval the recurrence coefficients $u_n(t)$ are bounded and
strictly positive $u_n(t)>0$. As is well known from general theory
of orthogonal polynomials \cite{Chi} the property $u_n>0$ for all
$n>0$ (together with reality of the coefficients $b_n$) guarantees
existence of a positive measure on the real axis. We already
constructed this measure explicitly \re{ort_cn}.

The remaining question is about uniqueness of the moment problem
for the case (i). Indeed, we have constructed explicitly the
orthogonality measure \re{ort_cn} corresponding to the moments \be
c_n(t) = \frac{d^n}{dt^n} \:
\left\{\frac{1}{\cn(wt\sqrt{e_1-e_3};k)}\right\}, \quad
n=0,1,2,\dots, \lab{moms_i} \ee where $t$ is assumed to belong to
the admissible interval \re{t_ad_cn}.

But in principle, it is possible that this measure is not unique.
Such situation is known as indeterminate moment problem
\cite{Akh_mom}, \cite{ShoTa}. In more details this means the
following. Assume that real moments $c_n$ are given and all
corresponding the Hankel determinants $D_n>0$ are positive. This
condition is equivalent to positivity of the recurrence
coefficients $u_n>0$ for $n=1,2,\dots$ and in turn, it guarantees
existence of a positive orthogonality measure on the real line
$-\infty < x <\infty$ (so-called the Hamburger moment problem
\cite{ShoTa}, \cite{ChiHamb}).  If this measure is unique (up to a
normalization condition) then the Hamburger moment problem is
called the {\it determinate}. If there exist at least two
different orthogonality measures then the Hamburger moment problem
is called indeterminate. In case of the indeterminate Hamburger
problem there exists infinitely many different measures (see
\cite{Akh_mom}, \cite{ShoTa} for details).

Finding criteria for determinacy of the Hamburger moment problem
is a nontrivial problem \cite{ChiHamb}. However, there is a simple
{\it sufficient} condition proposed by Carleman \cite{ShoTa}: if
\be \sum_{n=1}^{\infty} u_n^{-1/2} = \infty \lab{Carl} \ee then
the Hamburger problem is determinate. Of course, in \re{Carl} the
arithmetic value of the square root $\sqrt{u_n}$ is assumed.

We now show that the Hamburger problem for the moment problem
\re{moms_i} is determinate.

From \re{u_cn} it follows that for all $n$ and for fixed $t$ from
the admissible interval we have the inequalities
$$ 0<u_{2n}(t)<A(t) $$ where $$A(t)
= \frac{1}{\cn^2(u)} + \frac{k^2}{k'^2} = \frac{\dn^2(u)}{k'^2
\cn^2(u)}$$ is a fixed positive parameter (depending on $t$ but
not on $n$). Hence we have
$$
\sum_{n=1}^{\infty} u_{2n}^{-1/2} = \infty
$$
diverges. From the similar considerations it follows that
$$
\sum_{n=1}^{\infty} u_{2n+1}^{-1/2} = \infty
$$
Hence the Carleman condition \re{Carl} holds and we indeed have
the determinate moment problem. This means that the discrete
measure \re{ort_cn} is the only providing orthogonality of the
polynomials $P_n(z;t)$ on the real axis.

When $t$ approaches the endpoints of the admissible interval
\re{t_ad_cn} the recurrence coefficients $b_n(t), u_n(t)$ tend to
infinity. This "explosion" of the recurrence coefficients explains
corresponding "explosion" of the orthogonality measure when $t$
tends to he endpoints of the admissible interval.

Put now $t=0$ (the middlepoint of the admissible interval) . Then
it is seen (due to property $\eta_1+\eta_2+\eta_3=0$) that the
recurrence coefficient $b_n(0)$ vanishes \be b_n(0)= 0, \quad
n=0,1,2,\dots \lab{b_n_0} \ee and we arrive at a class of
so-called symmetric orthogonal polynomials with the recurrence
relation \cite{Chi} \be P_{n+1}(z) + u_n(0) P_{n-1}(z) = zP_n(z)
\lab{sym_rec} \ee For the recurrence coefficients $u_n(0)$ we have
\be u_{2n}(0) = 4 w^2 n^2 (e_1-e_2), \quad u_{2n+1}(0) = w^2
(2n+1)^2 (e_1-e_3). \lab{u_n_0} \ee These recurrence coefficients
correspond to the orthogonal polynomials introduced by Carlitz
\cite{Carlitz}. In turn, orthogonality relation for these OP
follows from the remarkable result by Stieltjes on presenting of
the Laplace transform of the Jacobi elliptic functions in terms of
continued fraction (for modern treatment of this and related
examples see e.g. \cite{Milne}. Today the orthogonal polynomials
introduced by Carlitz are called the Stieltjes-Carlitz orthogonal
polynomials related with elliptic functions \cite{Chi}. For
further details concerning these polynomials see \cite{LoBr}.

The orthogonality measure for polynomials corresponding to the
Stieltjes-Carlitz case \re{u_n_0} is obtained from our measure by
putting $t=0$. We thus have the orthogonality relation \be
\sum_{s=-\infty}^{\infty} {M_s(0) P_n(x_s;0)P_m(x_s;0)} = h_n(0)
\: \delta_{nm} \lab{ort_SC} \ee where the support of the measure
is the same, i.e. the points $x_s$ have the same expression
\re{loc_cn} and the concentrated masses are \be
M_s(0)=\frac{\pi}{k'K'} \: \frac{1}{v^{s-1/2} + v^{1/2-s}}
\lab{M_SC} \ee This measure was discovered by Stieltjes (Carlitz
showed that this measure provides orthogonality of the
corresponding polynomials $P_n(z)$). For the Stieltjes-Carlitz
polynomials the Hamburger moment problem is obviously determined
because $t=0$ belongs to the admissible interval.

We see that the Stieltjes-Carlitz polynomials appear naturally as
a very special case of the "elliptic Toda polynomials"
corresponding to "zero time" condition $t=0$.

We can finally summarize all these results as the
\begin{th}
Assume that $e_3<e_2<e_1$ are arbitrary real parameters with the
condition $e_1+e_2+e_3=0$. Assume that the recurrence coefficients
are given by formulas \re{b_cn}, \re{u_cn} with arbitrary positive
parameter $w$. Assume also that the parameter $t$ belongs to the
admissible interval \re{t_ad_cn}. Then the corresponding
orthogonal polynomials $P_n(x;t)$ are orthogonal on the uniform
grid \re{loc_cn} on the real axis with the concentrated masses
given by \re{M_cn}. The corresponding moments $c_n(t)$ are given
by \re{moms_i}. The moment problem is determinate.

\end{th}

\section{The orthogonality measure and recurrence coefficients for the case (ii)}
\setcounter{equation}{0} Consider now the case (ii). We have
analogously \be \tilde c_0(t) = \dc(iw\sqrt{e_1-e_3}t; k) =
\dn(w\sqrt{e_1-e_3}t; k') \lab{real_dn} \ee The Fourier series is
well known  \cite{WW} \be \dn(w \sqrt{e_1-e_3}t; k') = \frac{\pi
}{ K'} \: \sum_{n=-\infty}^{\infty} \frac{1}{v^{n} + v^{-n} }
\exp(\pi i n w t \sqrt{e_1-e_3}/K')  \lab{dn_F} \ee with the same
expression for $h$ as for the case (i). From considerations of the
first section we see that corresponding orthogonal polynomials
$P_n(z;t)$ have purely discrete measure located at the points \be
x_n = \frac{\pi w \sqrt{e_1-e_3}n}{K'}, \quad n=0, \pm 1, \pm 2 ,
\dots \lab{loc_dn} \ee Corresponding concentrated masses are \be
M_n(t)=\frac{2 \pi}{K'(v^n+v^{-n})} \: \exp(\pi w n
\sqrt{e_1-e_3}t/K') \lab{M_dn} \ee The admissible interval for $t$
is the same as for the case (i):
$$
-\frac{K}{w \sqrt{e_1-e_3}}< t < \frac{K}{w \sqrt{e_1-e_3}}
$$
As for the case (i) it is easily verified that inside the
admissible interval the recurrence coefficients $b_n(t),u_n(t)$
are real and $u_n >0$ which guarantees positivity of the measure:
$M_n(t) >0$ for all $n=0, \pm 1, \pm 2$ and for all $t$ from the
admissible interval.

The coefficients $u_n(t)$ has the expression \ba &&u_{2n}(t) = 4
w^2 n^2 (e_1-e_3) \: \left( k'^2 \: \frac{\sn^2(u)}{\cn^2(u)} +
\dn^2(2nu) \right)  \nonumber \\&& u_{2n+1}(t) = w^2 (2n+1)^2
(e_1-e_2) \: \left( \frac{1}{\cn^2(u)} + k^2 \:
\frac{\sn^2((2n+1)u)}{\dn^2((2n+1)u)} \right) \lab{udn} \ea The
recurrence coefficients $b_n(t)$ are expressed as \ba &&b_{2n}(t)
= w \sqrt{e_1-e_3} \:\left\{(2n+1) Z((2n+1) u +K) -\right.
\nonumber \\ &&  \left . 2n Z(2nu)   - (4n+1) (Z(u+K+iK') -
i\pi/(2K) ) \right\}  \nonumber \\ \lab{bdn} \\&&b_{2n+1}(t) = w
\sqrt{e_1-e_3} \:\left\{2(n+1) Z(2(n+1) u) -\right. \nonumber
\\ &&  \left . (2n+1) Z((2n+1)u +K)   - (4n+3) (Z(u+K+iK') - i\pi/(2K) )
\right\} \nonumber \ea

From \re{udn} it is clear that the coefficients $u_n(t)$ are
strictly positive $u_n(t)>0$ for any fixed value of the parameter
$t$ from the admissible interval.

From the same considerations it follows that the Hamburger moment
problem for the moments
$$
c_n(t) = \frac{d^n}{dt^n} \: \left\{
\frac{\dn(wt\sqrt{e_1-e_3}t;k)}{\cn(wt\sqrt{e_1-e_3}t;k)}
\right\}, \quad n=0,1,2,\dots
$$
is determinate for any value of the parameter $t$ from the
admissible interval \re{t_ad_cn}.

When $t=0$ then again, as in the case (i) the diagonal recurrence
coefficients are zero $b_n(0)=0$ and \be u_{2n} = 4w^2 n^2
(e_1-e_3), \quad u_{2n+1} = w^2 (2n+1)^2 (e_1-e_2) \lab{ub_dn_0}
\ee Note that the recurrence coefficients \re{ub_dn_0} are
obtained from the corresponding coefficients \re{u_n_0} of the
case (i) by a simple transposition $e_2 \leftrightarrow e_3$.
These recurrence coefficients \re{ub_dn_0} correspond to the
second class of the Stieltjes-Carlitz orthogonal polynomials
\cite{Carlitz}, \cite{LoBr} arising from the Laplace
transformation of the elliptic function $\dn(t;k')$.

We see, that just as in the case (i), the special choice of the
parameter $t=0$ leads to already known Stieltjes-Carlitz
orthogonal polynomials. For arbitrary values of the parameter $t$
from the admissible interval we obtain new orthogonal polynomials
with explicitly known recurrence coefficients and positive
discrete measure on a whole real axis.

We considered only two possible choices of the parameters $\beta,
\; q$ ($\beta=\omega_1/w, \; q= \omega_2$ and $\beta=\omega_1/w,
\; q= \omega_3$)  because for these two choices we get polynomials
$P_n(z;t)$ having positive orthogonality measure on the real axis.
It seems that all other values of the parameters $\beta,\: q$ (for
real values $w$) do not lead to polynomials with a positive
measure. Nevertheless, as we will see in the next sections, in the
degenerate cases of elliptic functions there are more
possibilities for these  parameters when the measure appears to be
a positive.

\section{Special choices of the parameter $t$}
\setcounter{equation}{0} The expressions \re{u_om} and \re{b_om}
show that for generic value of $t$ the recurrence coefficients
$b_n(t), u_n(t)$ are transcendent functions in $n$. We already saw
that for special choice $t=0$ we obtain $b_n(0)=0$ and $u_{2n}(0)$
as well $u_{2n+1}(0)$ are quadratic polynomials in $n$.

Put now \be t = \frac{M \: \omega_1}{N \: w}, \lab{om1_MN} \ee
where $M<N$ are mutually prime positive integers. Present the
number $n$ in the form
$$
n=Nr +s, \quad r=0,1,\dots, \; s=0,1,\dots, N-1
$$
Then for fixed $s=0,1,\dots, N-1$ the recurrence coefficients
$b_{2n}(t), \: b_{2n+1}(t)$ will be linear function in $r$  and
the coefficients $u_{2n}(t), \: u_{2n+1}(t)$ will be quadratic
polynomials in $r$ .

Indeed, from formulas \re{u_om}, \re{b_om}, using periodicity of
the Weierstrass functions, we get \ba &&u_{2Nj +
2s} = 4 w^2 (Nr+s)^2 \left(\epsilon_0 - \epsilon_1(s) \right), \\
\nonumber &&u_{2Nr+2s+1} = w^2 (2Nr + 2s +1)^2 \left(\epsilon_0 -
\epsilon_2(s)\right), \lab{sp_u} \ea where
$$
\epsilon_0= \wp((M+N)\omega_1/N), \quad \epsilon_1(s)=
\wp(\omega_j + 2sM\omega_1/N), \; \epsilon_2(s)= \wp(\omega_l +
(2s+1)M\omega_1/N)
$$
In these formulas $j=2, l=3$ for the case (i) and $j=3,l=2$ for
the case (ii).

Analogously, for the coefficients $b_{2n}(t)$ and $b_{2n+1}(t)$ we
obtain \ba &&b_{2Nr+2s}= w(-\eta_j + 2 \eta_1(Mr+Nr+s)) +w
(\kappa_1(s) (2Nr+2s+1)) - \nonumber \\ &&2 \kappa_2(s) (Nr+s) -
\kappa_0 (4(Nr+s)+1)) \lab{b1_sp} \ea \ba &&b_{2Nr+2s+1}=
w(-\eta_j
+ 2 \eta_1(Mr+3Nr+3s+2)) +w (\kappa_3(s) (2Nr+2s+2)) - \nonumber \\
&& \kappa_4(s) (2Nr+2s +1) - \kappa_0 (4(Nr+s)+3)) \lab{b2_sp} \ea
where
$$
\kappa_0=\zeta((M+N)\omega_1/N), \; \kappa_1(s) = \zeta(-\omega_l
+(2s+1) M\omega_1/N), \; \kappa_2(s) = \zeta(\omega_j +2s
M\omega_1/N)
$$
$$
\kappa_3(s) = \zeta(\omega_j +2s M\omega_1/N), \;
\kappa_4(s)=\zeta(-\omega_l +s M\omega_1/N)
$$
We see that indeed $u_{2Nr+2s}, u_{2Nr+2s+1}$ are quadratic in $r$
and $b_{2Nr+2s}, b_{2Nr+2s+1}$ are linear in $n$.

The constants $\epsilon_i(s), \kappa_i(s)$ can be found in a less
or more "explicit" form only for several values of $N,M$. We
already considered the case $t=0$ which corresponds to $M=0,\:
N=1, s=0$. Another simplest case corresponds to the choice
$$
M=1, N=2
$$
i.e. we choose $t = \frac{\omega_1}{2w}$. In this case the
parameter $s$ can take only 2 values $s=0,1$.

The zero moment becomes now $$c_0(t)=1/\cn(K/2;k) =
\sqrt{\frac{1+k'}{k'}}.$$

Using relations  and \re{zeta_half} we can calculate the
recurrence coefficients $b_n(\omega_1/(2w))$ \be
b_{4n}(\omega_1/(2w))=  w \: \sqrt{e_1-e_3} (2n (k'+3) +1)
\lab{b_4n} \ee

\be b_{4n+1}(\omega_1/(2w))= w \: \sqrt{e_1-e_3} (2n (3k'+1)
+2k'+1) \lab{b_4n+1} \ee

\be b_{4n+2}(\omega_1/(2w))= w \: \sqrt{e_1-e_3} (2n (3k'+1) +4k'
+1) \lab{b_4n+2} \ee

\be b_{4n+3}(\omega_1/(2w))= w \: \sqrt{e_1-e_3} (2n (k'+3) + 2k'
+5) \lab{b_4n+3} \ee

Similarly, using relations \re{wp_half} we obtain expressions for
the coefficients $u_n(\omega_1/(2w))$;  \be
u_{2n+1}(\omega_1/(2w)) = w^2 \: (e_1-e_3) (2n+1)^2 \: 2k'
\lab{u_odd} \ee

\be u_{4n}(\omega_1/(2w)) = 16 w^2 (e_1-e_3) n^2 k'(1+k')
\lab{u_4n} \ee

\be u_{4n+2}(\omega_1/(2w)) = 4 w^2 (e_1-e_3) (2n+1)^2 (1+k')
\lab{u_4n+2} \ee

Note that the combination $w \sqrt{e_1-e_3}$ plays the role of a
scaling parameter, hence we can put $w=1, \; e_1-e_3 =1$ without
loss of generality.

We then have
\begin{th}
Assume that $0<k'<1$ is an arbitrary parameter.  Let the
recurrence coefficients  $b_n$ be defined as \ba &&(k'+3)n/2 + 1,
\quad \mbox{if} \quad n=0 \quad (mod \quad 4); \nonumber
\\&&(3k'+1)n/2 + (k'+1)/2, \quad \mbox{if} \quad n=1 \quad (mod \quad 4);
\nonumber \\ &&(3k'+1)n/2+k', \quad \mbox{if} \quad n=2 \quad (mod
\quad 4); \nonumber   \\&&(k'+3)n/2 +(k'+1)/2 \quad \mbox{if}
\quad n=3 \quad (mod \quad 4) \lab{b_spec} \ea and the recurrence
coefficients $u_n$ be defined as \ba && k'(k'+1) n^2, \quad
\mbox{if} \quad n=0 \quad (mod \quad 4); \nonumber \\&& 2k' n^2,
\quad \mbox{if} \quad n=1,3 \quad (mod \quad 4); \nonumber
\\&&(1+k')n^2, \quad \mbox{if} \quad n=2 \quad (mod \quad 4)
\lab{u_spec} \ea

Then corresponding monic orthogonal polynomials $P_n(x)$ are
orthogonal with purely discrete measure on the real line \be
\sum_{s=-\infty}^{\infty} M_s P_n(x_s) P_m(x_s) = h_n \:
\delta_{nm}, \lab{spec_ort} \ee where the grid of orthogonality is
$$x_s = \frac{\pi (s-1/2)}{K'}$$ and corresponding concentrated
masses are
$$
M_s = \frac{\pi}{k'K'} \: \frac{v^{(s-1/2)/2}}{v^{s-1/2} +
v^{1/2-s}}, \quad v = \exp(-\pi K/K')
$$
The normalization constants $h_n$ are $$h_n = c_0 u_1 u_2 \dots
u_n$$ where
$$
c_0 =\sqrt{\frac{1+k'}{k'}}
$$

\end{th}

\section{Degenerated cases}
\setcounter{equation}{0} Consider degenerated cases of obtained
orthogonal polynomials. These degenerated cases arise when two or
three of the parameters $e_i$ coincide. Using our choice
$e_1>e_2>e_3$ we see that there are two possibilities when two of
the parameters coincide:

(i) $e_1=e_2=a, \; e_3=-2a$, where $a$ is a positive parameter. In
this case the real period $2 \omega_1$ tends to infinity, whereas
the imaginary period remains finite
$$
\omega_3 = \frac{\pi i}{\sqrt{12 a}}
$$
Without loss of generality we can take $a=1/3$ (changing of $a$
leads only to scaling of the argument $z$ of corresponding
functions). The Weierstrass functions are then reduced to
hyperbolic ones, e.g.
$$
\wp(z) \to 1/3 +\frac{1}{\sinh^2(z)}
$$
The modular parameter becomes $k=1$ and the Jacobi elliptic
functions become hyperbolic as well: $\sn(z;1) = \tanh(z), \;
\cn(z;1) = \dn(z;1)=1/\cosh(z)$. The imaginary period in this case
is $2 \omega_3 = i \pi$.

(ii) $e_2=e_3=-a, \; e_1 =2a$ with some positive parameter $a$.
Then the imaginary period $2 \omega_3$ becomes infinity whereas
the real period is finite
$$
\omega_1 = \frac{\pi}{\sqrt{12a}}
$$
Again we can put $a=1/3$, then the Weierstrass  function $\wp(z)$
is degenerated to trigonometric form:
$$
\wp(z) \to -1/3 + \frac{1}{\sin^2(z)}
$$
The modular parameter $k=0$ in this limit and we have
$\sn(z;0)=\sin(z), \; \cn(z;0) = \cos(z), \; \dn(z;0) = 1 $.

Consider first the degenerated cases (i) and (ii) of the elliptic
polynomials obtained in the previous section. In the "hyperbolic"
limit $k=1$ the functions $c_0(t)$ becomes the same $c_0(t) =
\cosh(t)$. This case is non-interesting because it corresponds to
degeneration of the orthogonal polynomials: the Hankel
determinants become zero $D_n(t)$ for infinitely many $n$. In
turn, this corresponds to zero recurrence coefficients for even
$2n$: $u_{2n}(t)=0$ as can be seen from explicit formulas for
$u_n(t)$ in the hyperbolic limit.

In the trigonometric limit $k=0$ again both functions $c_0(t)$
coincide: $c_0(t)=1/\cos(t)$. This function $c_0(t)$ corresponds
to some elementary  solutions of the restricted Toda chain;
corresponding orthogonal polynomials $P_n(z;t)$ coincide with a
special case of the Meixner-Pollaczeck polynomials.

Indeed, for this case it is elementary verified that the
recurrence coefficients have the expression: \be u_n(t) =
\frac{n^2}{\cos^2 t}, \quad b_n(t)= (2n+1) \: \tan t \lab{MPT_ub}
\ee On the other hand the monic Meixner-Pollaczeck polynomials
\cite{KS} \be P_n^{(\lambda)}(x;\phi) =\frac{(2 \lambda)_n
e^{in\phi} }{(2\sin \phi)^n}{_2}F_1\left( {-n, \lambda+ix \atop 2
\lambda}; 1-e^{-2i\phi} \right) \lab{st_MP} \ee (here
$(a)_n=a(a+1)\dots (a+n-1)$ is the standard shifted factorial)
depend on two parameters $\lambda, \phi$ and have the recurrence
coefficients \be u_n =\frac{n(n+2 \lambda-1)}{4 \sin^2\phi}, \quad
b_n = - \frac{n+ \lambda}{\tan \phi} \lab{ub_MP} \ee Comparing the
recurrence coefficients we see that our polynomials coincide with
(rescaled) Meixner-Pollaczeck polynomials with the parameters
$\lambda= 1/2, \; \phi= t + \pi/2$.

The explicit expression for our polynomials $P_n(x;t)$ is
 \be P_n(x;t)= \frac{n! i^n
e^{int} }{\cos^nt} \;  {_2}F_1\left( {-n, 1/2+ix/2 \atop  1};
1+e^{-2it} \right) \lab{our_MP} \ee From the standard formulas for
the Meixner-Pollaczeck polynomials \cite{KS} we find that that our
polynomials $P_n(x;t)$ are orthogonal on whole real axis: \be
\int_{x=-\infty}^{\infty} P_n(x;t) P_m(x;t) W(x;t) dx = h_n(t) \:
\delta_{nm}, \lab{ort_MP} \ee where the weight function $W(x,t)$
is \be W(x,t) =\frac{\pi \: e^{tx}}{2\cosh(\pi x/2)} \lab{W_MP}
\ee The weight is well defined provided that $t$ belongs to the
admissible interval $-\pi/2 < t < \pi/2$. When $t$ is inside this
interval all moments $c_n(t)$ exist and we have the normalization
condition
$$
\int_{-\infty}^{\infty} W(x;t) dx =  \int_{-\infty}^{\infty}
\frac{\pi \: e^{tx}}{2\cosh(\pi x/2)} dx= c_0(t) = \frac{1}{\cos
t}
$$
It is instructive to see how the continuous orthogonality relation
\re{ort_MP} arises in the limiting case of the discrete-type
orthogonality relation \re{ort_cn}. Indeed, in the trigonometric
limit $k \to 0$ we have $K(k) \to \pi/2$ and $K'(k) \to \infty$.
So (recall that we assume $e_1=2/3, e_2=e_3=-1/3$) from
\re{loc_cn} we see that the grid step $\Delta x(s) =x_{s+1}-x_s $
becomes infinitely small and the sum in rhs of \re{ort_cn} becomes
an integral \re{ort_MP} after appropriate definition of the
continuous variable $x$.

Note also that the special case considered in the previous section
(recurrence coefficients given by \re{b_spec}, \re{u_spec}) in the
limit $k=0$ corresponds to the formulas \re{MPT_ub} for $t=\pi/4$:
$$
u_n = 2n^2, \quad b_n = 2n+1.
$$

Consider now more general class of degenerated solutions
corresponding to the case when $\beta=0, \mu_1=wq/3, \mu_0=0$ and
$q$ is an arbitrary real parameter. We then have in the hyperbolic
limit ($k=1$) \be c_0(t)= \frac{\sinh(wt+q)}{\sinh(q) \sinh(wt)} =
\frac{1}{2(e^{2q}-1)} + \frac{2}{1-e^{-2wt}}  \lab{c_cth} \ee
Introduce the function \be c_0^{(0)}(t) = \frac{2}{1-e^{-2wt}} =
\frac{e^{wt}}{\sinh(wt)}.\lab{c00_M} \ee The function \re{c_cth}
differs from $c_0^{(0)}(t)$ only by adding of a term
$e^{-q}/\sinh(q)$ not depending on $t$. We have
$$
c_n(t) = c_n^{(0)}(t)= \frac{d^n c_0^{(0)}(t)}{dt^n}, \quad
n=1,2,\dots
$$
Thus all moments corresponding to functions $c_0^{(0)}(t)$ and
$c_0(t)$ coincide apart from zero moments.

Introduce linear functional $\sigma(t)$ and $\sigma^{(0)}(t)$ by
their moments
$$
\langle \sigma(t), x^n \rangle = c_n(t), \quad \langle
\sigma^{(0)}(t), x^n \rangle = c_n^{(0)}(t), \quad n=0,1,2,\dots
$$
We see that the functionals $\sigma(t)$ and $\sigma^{(0)}(t)$ are
related as \be \sigma^{(0)}(t) = \sigma(t) + e^{-q}/\sinh(q)
\delta_0, \lab{ss_delta} \ee where $\delta_0$ is the Dirac
delta-functional, corresponding to inserting a unit concentrated
mass to the point $x=0$:
$$
\langle \delta_0, x^n \rangle = \delta_{n0} $$ Assume that
$\rho^{(0)}(x)$ is the orthogonality measure for the polynomials
$P_n^{(0)}(z;t)$ corresponding to the functional
$\sigma^{(0)}(t)$:
$$
\int_{-\infty}^{\infty} P_n^{(0)}(x;t) P_m^{(0)}(x;t) d
\rho^{(0)}(x;t) = h_n^{(0)} \: \delta_{nm}
$$
Then the orthogonality measure corresponding to he polynomials
$P_n(x;t)$ is
$$
\rho(x;t) = \rho^{(0)}(x;t) + e^{-q}/\sinh(q) \: \delta(x)
$$
Thus indeed the weight of orthogonality for the polynomials
$P_n(x;t)$ is obtained from the corresponding orthogonality weight
for the polynomials $P_n^{(0)}(x;t)$ by adding of a concentrated
mass $\coth (q)$ at the point $x=0$.

Now we show that the function $c_0^{(0)}(t)$ given by \re{c00_M}
generates a special class of the Meixner polynomials. Indeed, it
is elementary verified that two sequences \be b_n(t) =
-\frac{2w(n+1 + ne^{2wt})}{e^{2wt}-1}, \quad u_n = \frac{4 w^2 n^2
e^{2wt}}{(e^{2wt}-1)^2}, \quad n=0,1,2,\dots \lab{bu_M} \ee
satisfy the restricted Toda chain equations \re{Toda} together
with the condition $b_0(t) = \dot c_0^{(0)}(t)/c_0^{(0)}(t)$. Thus
the recurrence coefficients \re{bu_M} correspond to orthogonal
polynomials $P_n^{(0)}(z;t)$ having the moments $c_n^{(0)}(t) =
d^n c_0^{(0)}(t)/dt^n$. On the other hand, we can easily identify the recurrence coefficients
\re{bu_M} with the those for the special class of the Meixner polynomials.

Indeed, the Meixner polynomials $P_n(x;\beta;c)$ have two real
parameters $\beta,c$ and have the recurrence coefficients
\cite{KS} \be b_n = \frac{n+(n+\beta)c}{1-c}, \quad u_n =
\frac{c}{(1-c)^2}n(n+\beta-1) \lab{bu_Me} \ee Explicitly the
Meixner polynomials are expressed in terms of the Gauss
hypergeometric function \cite{KS} \be P_n(x;\beta;c) =\kappa_n \:
{_2}F_1\left(  {-n, -x  \atop  \beta};  1-1/c \right) \lab{Me_P}
\ee where $\kappa_n$ is an appropriate normalization factor to
provide monicity of the polynomials $P_n(x;\beta;c)$. The Meixner
polynomials are orthogonal on the uniform semi-infinite grid
\cite{KS}: \be \sum_{s=0}^{\infty} \frac{c^s (\beta)_k}{k!}
P_n(k;\beta;c) P_m(k;\beta;c) = h_n \: \delta_{nm} \lab{ort_Me}
\ee Obviously we should have $0<c<1$ to provide positivity
property of the measure.

Comparing the recurrence coefficients \re{bu_M} with \re{bu_Me} we
see that $\beta=1, \; c=e^{-2wt}$ and polynomials $P_n^{(0)}(x;t)$
coincide with the corresponding rescaled Meixner polynomials: \be
P_n^{(0)}(x;t) = (-2w)^n \: P_n(-x/(2w); 1; e^{-2wt}) \lab{P_PM}
\ee Orthogonality relation for the polynomials $P_n^{(0)}(x;t)$
looks as \be \sum_{s=0}^{\infty} e^{-2swt}  P_n^{(0)}(-2ws;t)
P_m^{(0)}(-2ws;t) = h_n^{(0)}(t) \: \delta_{nm} \lab{ort_P0} \ee
To provide positivity of the measure for $t>0$ we should have
$w>0$. Thus polynomials $P_n^{(0)}(x;t)$ are orthogonal on the
uniform grid of the negative real axis.

Return to the polynomials $P_n(z;t)$ corresponding to the function \re{c00_M}.
The recurrence coefficients for the polynomials $P_n(z;t)$ are obtained from the recurrence coefficients \re{b_FS},
\re{u_FS} by the limiting procedure $e_2 \to e_1$:
\be
b_n(t) = w(n+1) \coth(w(n+1)t+q) - wn \coth(wnt + q) - w(2n+1) \coth(wt) \lab{b_M} \ee
and
\be
u_n(t) = w^2n^2 \: \frac{\sinh(w(n+1)t + q) \sinh(w(n-1)t + q)}{\sinh^2(wnt + q) \sinh^2(wt)} \lab{u_M} \ee

As we already showed, the orthogonality relation for the
polynomials $P_n(x;t)$ corresponding to the function \re{c00_M} is
obtained from \re{ort_P0} by adding of a concentrated mass at the
point $x=0$. Explicitly we have \be \sum_{s=0}^{\infty} 2
e^{-2wst} P_n(-2ws,t)P_m(-2ws,t) + M P_n(0;t)P_m(0;t) = h_n \:
\delta_{nm} \lab{ort_M} \ee where the value of the mass inserted
at $x=0$ is $$M=\frac{e^{-q}}{ \sinh(q)}.$$ It is assumed that
$w,t>0$ in order to provide convergence of series in lhs of
\re{ort_M}.

The normalization constants $h_n(t)$ are expressed through the
recurrence coefficients \re{u_M} as \be h_n(t) = c_0(t) u_1(t)
u_2(t) \dots u_n(t) \lab{h_M} \ee Note that for $n=m=0$ formula
\re{ort_M} gives an obvious identity
$$
\sum_{s=0}^{\infty} e^{-2wts} + \frac{e^{-q}}{ \sinh(q)} = \frac{\sinh(wt+q)}{\sinh(q) \sinh(wt)} = c_0(t)
$$

When $q \to \infty$ we see that $M \to 0$ and  $c_0(t) \to
\frac{e^{wt}}{\sinh(wt)} = c_0^{(0)}(t)$  hence in this limit the
polynomials $P_n(z;t)$ become the ordinary Meixner polynomials
\re{P_PM}.

The polynomials obtained by an adding of a concentrated mass at
the point $x=0$ of the orthogonality measure for the Meixner
polynomials are called the modified Meixner polynomials and were
proposed by R.Askey as an interesting object for further
investigations \cite{Askey}. Properties of these polynomials were
intensively studied in \cite{AMA} and \cite{Bavinck}. In
particular it was shown that these polynomials satisfy difference
equations of finite and infinite order.

\section{Trigonometric limit}
\setcounter{equation}{0} Consider the trigonometric limit when
$e_1=2/3, e_2=e_3=-1/3$. Put $\beta=0, \mu_0=0, \mu_1 = - wq/3.$
The recurrence coefficients take the form \be b_n(t) = w(n+1)
\cot(w(n+1)t +q) - wn \cot(wnt+q) -(2n+1)w\cot(wt) \lab{b_trig}
\ee and \be u_n(t) = w^2 n^2 \: \frac{\sin((n+1)wt +
q)\sin((n-1)wt + q)}{\sin^2(nwt + q) \sin^2(wt)} \lab{u_trig} \ee
and the function $c_0(t)$ is \be c_0(t) =
\frac{\sin(wt+q)}{\sin(q) \sin(wt)} \lab{c0_trig} \ee From the
expression \re{u_trig} we see that for real value $t,w$ the
coefficient $u_n(t)$ cannot have the same sign for all $n$. Hence,
corresponding measure for orthogonal polynomials is not positive
definite.

Nevertheless, there is one interesting special case leading to a positive definite measure on a
{\it finite set} of points on the real axis. Indeed, put $w=1, q=\pi/2$, then
\be
c_0(t) = \cot(t) \lab{c0_cot} \ee
and
$$
b_n(t) = -(n+1) \tan((n+1)t) + n \tan(nt) - (2n+1) \cot(t),$$
$$u_n(t) = n^2 \: \frac{\cos((n+1)t) \cos((n-1)t)}{\cos^2(nt) \sin^2(t)}
$$
Assume that
$$
t=\tau=\frac{\pi}{2(N+2)}
$$
 for some positive integer $N=2,3,\dots$. Then it is seen that $u_k>0$ for $1\le k \le N$ and $u_{N+1}=0$.
 This condition guarantees that the finite set of polynomials $P_0(x;\tau), P_1(x;\tau), \dots, P_{N}(x;\tau)$
 will be orthogonal on the set of (real) zeros $x_s$ of the polynomial $P_{N+1}(x;t)$:
 \be
 \sum_{s=0}^N \rho_s P_n(x_s;\tau) P_m(x_s; \tau)= h_n \: \delta_{nm} \lab{fin_ort} \ee
with postive weights $\rho_s$.

However explicit expressions for zeros $x_s$ and the weights $\rho_s$ in this case are still unknown.

\section{Completely degenerated case. The Krall-Laguerre polynomials}
\setcounter{equation}{0}
Finally, consider the case when all roots coincide $e_1=e_2=e_3=0$.
Then the elliptic functions are degenerated to simple rational ones:
$\wp(z) = 1/z^2, \: \zeta(z) =1/z, \: \sigma(z) =z$.

Assuming again $w=1, \beta=0$ we have
\be
c_0(t) = \frac{t+q}{qt} = 1/t + 1/q \lab{c_0_rat} \ee
The first term $1/t$ in rhs of \re{c_0_rat} generates the Laguerre polynomials. Indeed,
it is elementary verified that $$
b_n(t) = -\frac{2n+1}{t}, \quad u_n(t)=\frac{n^2}{t^2}
$$
is a solution of the restricted Toda chain corresponding to the
initial condition $c_0^{(0)}(t) = 1/t$. These recurrence
coefficients correspond to the Laguerre polynomials
$L_n^{(0)}(-xt)$ \cite{KS}.

The second constant term in \re{c_0_rat} describes adding of a
concentrated mass to the measure at the endpoint $x=0$ of the
orthogonality interval. We thus obtain that the polynomials
$P_n(z;t)$ corresponding to \re{c_0_rat} coincide with the
so-called Krall-Laguerre polynomials (see, e.g. \cite{Krall} for
details). The measure for the latter is obtained by the
adding of (an arbitrary) concentrated mass to the point $x=0$ of
the orthogonality interval for the Laguerre polynomials
$L_n^{(0)}(x)$. In our case the orthogonality relation looks as
\be \int_{-\infty}^{0} P_n(x;t) P_m(x;t) e^{xt} dx +
\frac{P_n(0;t) P_m(0;t)}{q} = h_n(t) \: \delta_{nm}
\lab{ort_Lag_q} \ee The recurrence coefficients are \ba &&b_n(t) =
-\frac{n}{nt+q} + \frac{n+1}{(n+1)t +q} - \frac{2n+1}{t} \nonumber
\\ &&u_n(t) = \frac{n^2}{t^2} \: \frac{((n-1)t +q)((n+1)t
+q)}{(nt+q)^2} \lab{bu_Lag_q} \ea

Remarkably enough that the Krall-Laguerre polynomials belong to a
class of 3 families of "non-classical" orthogonal polynomials
satisfying the ordinary eigenvalue problem for the linear
differential operator of the 4-th order \cite{Krall}.

\section{Continued fractions and the Hankel determinants connected
with the Jacobi elliptic functions} \setcounter{equation}{0} As we
saw, any solution $u_n(t), b_n(t)$ of the restricted Toda chain is
generated by the only function $c_0(t)$.

The Stieltjes function $F(z;t)$ for the corresponding orthogonal
polynomials $P_n(z;t)$ is given by the Laplace transform \be
F(z;t) = \frac{1}{c_0(t)} \: \int_0^{\infty} c_0(y+t) \: e^{-yz}
dy \lab{Lap_yz} \ee (The factor $1/c_0(t)$ in front of the
integral in \re{Lap_yz} is needed to provide "conventional"
asymptotic behavior $F(z;t) = z^{-1} + O(z^{-2})$). On the other
hand, it is well known \cite{Chi} that any Stieltjes function with
such asymptotic behavior generates a continued fraction of the
Jacobi type (so-called J-type continued fraction): \be
 F(z) = {1\over\displaystyle z-b_0 - {\strut u_1
\over\displaystyle z-b_1 - {\strut u_2 \over {z-b_2  -  \dots
}}}}, \qquad  \lab{cont_Jac} \ee where $b_n,u_n$ are corresponding
recurrence coefficients for the polynomials $P_n(z;t)$. Thus we
can construct explicitly families of continued fractions
\re{cont_Jac} starting from known solution $u_n(t), b_n(t)$ of the
Toda chain corresponding to the function $c_0(t)$ given by
\re{c0_FS}.

For a special choice of the parameters $\beta, q, \mu_1$ we can
obtain families of the Jacobi elliptic functions $\sn(t),
\cn(t),\dn(t)$ as well as related functions obtained by simplest
modular transforms. Thus we can generalize the Stieltjes results
(extended and generalized by Milne \cite{Milne}) who obtained
continued fractions corresponding to the Laplace transform of the
Jacobi elliptic functions in a special case $t=0$ of the formula
\re{Lap_yz}.

Moreover, from our results it follows a simple formula for the
corresponding Hankel determinants $D_n(t)$ defined by \re{DeltaT}.

Indeed, by \re{uh} we have \be D_n(t) = h_0(t) h_1(t) \dots
h_{n-1}, \quad n=1,2,3,\dots \lab{D_h} \ee or, equivalently, \be
D_n(t) = c_0(t) u_1^{n-1}(t) u_2^{n-2}(t) \dots u_{n-2}^2(t)
u_{n-1}(t) \lab{D_u} \ee Using explicit formulas \re{b_FS} and
\re{u_FS} we find \be h_n(t) = e^{\mu_1(t+\beta) + \mu_0} n!^2
w^{2n} \: \frac{\sigma((n+1) w(t+\beta) +q)}{\sigma(n w(t+\beta)
+q) \sigma^{2n+1}(wt+\beta)}. \lab{h_FS} \ee   Note that the
parameter $\mu_0$ is inessential, because it doesn't contribute to
the recurrence coefficients $b_n(t), u_n(t)$, nevertheless such
parameter is convenient when we would like to take $c_0(t)$
coinciding with prescribed Jacobi functions. For the Hankel
determinants we have the expression  \be D_n(t) = \kappa_n \:
\frac{\sigma(wn(t+\beta) + q)}{\sigma(q) \sigma^{n^2}(w(t+\beta))
} \: \exp((\mu_1(t+\beta) + \mu_0)n) \lab{D_n_FS} \ee where \be
\kappa_n = 1!^2 2!^2 \dots (n-1)!^2 w^{n(n-1)} \lab{kap_n} \ee

Consider now 3 special cases corresponding to the basic Jacobi
elliptic functions $\sn(t), \cn(t), \dn(t)$. In all these case we
can assume that $e_1-e_3=1$. Indeed, the Jacobi elliptic functions
$\sn(t,k), \cn(t,k), \dn(t,k)$ depend on the modulus
$$
k^2 = \frac{e_2-e_3}{e_1-e_3}.
$$
Hence we can pass from given parameters $e_i$ to the scaling
parameters $\gamma \: e_i, \: i=1,2,3$ with some nonzero constnat
$\gamma$. Such transformation doesn't change the Jacobi elliptic
functions. Hence we can always can assume that $e_1-e_3=1$.

For the function $\sn(t)$ we put $w=1, \beta=q =\omega_3$. We then
have the recurrence coefficients \ba &&
u_{2n}(t) = 4n^2 k^2 \left( \sn^2(t) - \sn^2(2nt) \right)\nonumber \\
&&u_{2n+1}(t) = (2n+1)^2 \left( k^2\sn^2(t) -
\frac{1}{\sn^2((2n+1)t)} \right) \nonumber \ea and \ba &&b_{2n}(t)
= (2n+1) \left\{ Z((2n+1)t + iK') + i\pi/(2K)\right\} - 2n Z(2nt)
-(4n+1) Z(t)    \nonumber \ea \ba &&b_{2n+1}(t) = (2n+2)
Z((2n+2)t,k )  - (2n+1) (Z((2n+1)t + iK',k) + i\pi/(2K)) -
\nonumber \\&&(4n+3) Z(t,k) \nonumber \ea

For the function $\cn(t)$ we put $w=1, \:  \beta=\omega_3,\; q
=\omega_2$. We then have the recurrence coefficients \ba &&
u_{2n}(t) = 4n^2 k^2 \left( -\cn^2(t) + k'^2 \: \frac{ \sn^2(2nt)}{\dn^2(2nt)} \right)\nonumber \\
&&u_{2n+1}(t) = (2n+1)^2 \left( -\dn^2(t) - k'^2 \:
\frac{\sn^2((2n+1)t)}{\cn^2((2n+1)t)} \right) \nonumber \ea and
\ba &&b_{2n}(t) = (2n+1) \left\{ Z((2n+1)t,k) -
\frac{\sn((2n+1)t,k) \dn((2n+1)t,k) }{\cn((2n+1)t,k)} \right\}-
\nonumber \\&&  2n Z(2nt + K,k) -(4n+1) Z(t,k)  \nonumber \ea

\ba &&b_{2n+1}(t) = (2n+2) Z((2n+2)t,k) +K) + \nonumber \\&&(2n+1)
\left\{-Z(2n+1)t,k)  + \frac{\sn((2n+1)t,k) \dn((2n+1)t,k)
}{\cn((2n+1)t,k)} \right \}  - \nonumber
\\&& -(4n+3) Z(t,k)  \nonumber \ea

For the function $\dn(t)$ we put $w=1, \: \beta=\omega_3,\; q
=\omega_1$. We then have the recurrence coefficients \ba &&
u_{2n}(t) = 4n^2 \left( -\dn^2(t) - k'^2 \: \frac{ \sn^2(2nt)}{\cn^2(2nt)} \right)\nonumber \\
&&u_{2n+1}(t) = (2n+1)^2  k^2 \left( -\cn^2(t) + k'^2 \:
\frac{\sn^2((2n+1)t)}{\dn^2((2n+1)t)} \right) \nonumber \ea and
\ba &&b_{2n}(t) = (2n+1) Z((2n+1)t + K,k) -   2n \left\{ Z(2nt,k)
- \frac{\sn(2nt) \dn(2nt)}{\cn(2nt)}  \right\}
 - (4n+1) Z(t,k) \nonumber \ea

\ba &&b_{2n+1}(t) = (2n+2) \left\{ Z((2n+2)t,k) -
\frac{\sn((2n+2)t) \dn((2n+2)t)}{\cn((2n+2)t)}      \right\} -
\nonumber \\&& (2n+1) Z(2nt+K,k)
 - (4n+3) Z(t,k) \nonumber \ea

\section{Concluding remarks}
\setcounter{equation}{0}
The function $c_0(t)$ given by \re{c0_FS} is closely related with simplest solutions
of the Lam\'e equation.

Indeed, consider the Lam\'e equation in the form \cite{Akhiezer}
\be \frac{d^2 y}{du^2} = \left\{n(n+1) \: \wp(u) + l  \right\} y
\lab{Lame} \ee In case if $n$ is a positive integer one can
construct explicit solutions of the Lame equation in the form
\cite{Akhiezer} \be \phi(u) = e^{\lambda u} \: \frac{\sigma(u-a_1)
\dots \sigma(u-a_n)}{\sigma^n(u)} \lab{Lame_n_sol} \ee where the
constants $\lambda, a_1, \dots, a_n$ can be determined from the
Lam\'e equation. In the simplest nontrivial case $n=1$ we have \be
\phi(u) =\frac{\sigma(u+a)}{\sigma(u)} \: \exp(-\zeta(a) u)
\lab{Lame_1_sol} \ee where the parameter $a$ is related with the
spectral parameter $l$ of the Lam\'e equation by the
transcendental equation $\wp(a)=l$. In case if $a \ne \omega_k, \:
k=1,2,3$ we have the second linearly independent solution of the
Lam\'e equation in the form \be \phi(u)
=\frac{\sigma(u-a)}{\sigma(u)} \: \exp(\zeta(a) u)
\lab{Lame_1_sol2} \ee We see that our function $c_0(t)$ \re{c0_FS}
coincides with the Lam\'e solution \re{Lame_1_sol} when $\mu_1 =
-\zeta(a)$.

This means that if the Stieltjes function $F(z;t)$ is the Laplace
transform of the solution of the Lam\'e equation \be F(z;t) =
\int_0^{\infty} \left\{ \frac{\sigma(u+a+t)}{\sigma(u+t)} \:
e^{-(\zeta(a) +z)u} \right \} du \lab{St_Lame} \ee then
corresponding orthogonal polynomials $P_n(z;t)$ will have the
recurrence coefficients given by \re{b_FS} and \re{u_FS} where
$w=1, \: \beta=0, q =a, \mu_1 = -\zeta(a)$.

In \cite{Chu1}, \cite{Chu2} some explicit continued fractions
connected with the Lam\'e solutions \re{Lame_1_sol} and
\re{Lame_n_sol} were announced without any proof or even idea of
proof. The authors of \cite{Chu1}, \cite{Chu2} considered
Stieltjes functions of kind of \re{St_Lame} but the argument of
these function was $w=\wp(a)$ instead of $z$. This leads to
explicit continued fractions which do not resemble presented in
the present paper. It would be desirable to connect results in
\cite{Chu1} and \cite{Chu2} with our ones.

Another possible generalization consists in passing to the
associated polynomials. Indeed, we considered here only solutions
for the restricted Toda chain, i.e. under the condition $u_0=0$.
Nevertheless, solutions \re{b_FS}, \re{u_FS} can be easily
extended to the non-restricted case if one replace $n$ with $n+c$
in corresponding formulas, where $c$ is an arbitrary constant not
depending on $t$. Then we obtain solution of the nonrestricted
Toda chain if $c \ne 0, \pm 1, \pm 2, \dots$. Such replacement
$b_n \to b_{n+c}, \; u_n \to u_{n+c}$ is well known and leads to
replacing of the orthogonal polynomials $P_n(z;t)$ with their {\it
c-associated} polynomials. Valent already considered
\cite{Valent3}, \cite{Valent1} the $c$-associated polynomials
corresponding to the Stieltjes-Carlitz polynomials. He was able to
find an explicit orthogonality measure in some special cases. In
our case (i.e. for $t \ne 0$) the corresponding analysis seems to
be much more complicated.

Nevertheless, there is a simple special case of the associated
polynomials for $t\ne 0$ which leads again to solutions of the
restricted Toda chain.

Indeed, take $c=1$ in formulas \re{u_FS} and \re{b_FS} and then
put $q=0, \: w=1, \beta=0, \mu_1=0$. We obtain the recurrence
coefficients \be b_n(t) = (n+2) \zeta((n+2)t) -(n+1) \zeta((n+1)t)
- (2n+3) \zeta(t) \lab{b_W} \ee and \be u_n(t)=(n+1)^2 \:
\left(\wp(t) - \wp((n+1)t)   \right) = \frac{\sigma(nt)
\sigma((n+2)t)}{\sigma^2(t) \sigma^2((n+1)t)} \lab{u_W} \ee It is
seen from \re{u_W} that $u_0(t)=0$, hence we deal again with a
solution of the restricted Toda chain. From \re{b_W} we have
$$
b_0(t) = 2 \zeta(2t) - 4\zeta(t) = \dot c_0(t)/c_0(t)
$$
whence \be c_0(t) = \dot \wp(t), \lab{dot_w} \ee where we used the
identity \cite{WW}
$$
\zeta(2z) = \zeta(z) + \frac{\wp''(z)}{2 \wp'(z)}
$$
Thus we obtained that the function $c_0(t) = \dot  \wp(t)$
generates another solution of the restricted Toda chain described
by formulas \re{b_W}, \re{u_W}. This solution was already
presented by Chudnovsky brothers in \cite{Chu} The corresponding
orthogonal polynomials $P_n(z;t)$ seems not to possess positivity
property for the Hankel determinants $D_n(t)$, hence their measure
will not be positive on the real axis.

Nevertheless, if $t=2 \omega_1/(N+2)$ with some positive integer
$N$ we have $u_{N}(t)=0$. This means that  polynomials $P_n(z;t)$
will be orthogonal on a finite set of points $x_s$: \be
\sum_{s=0}^{N-1} w_s P_n(x_s;t) P_m(x_s;t) = h_n(t) \:
\delta_{nm}, \quad n,m =0,1,\dots, N-1 \lab{fin_W} \ee where $x_s$
are roots of the polynomial $P_N(x)$:
$$
P_N(x_s) =0, \quad s=0,1,\dots, N-1
$$
Assume that all the roots $x_s$ are simple. Then the discrete
weight function $w_s$ can be presented in the form \cite{Chi} \be
w_s = \frac{h_{N-1}}{P_{N-1}(x_s) P'_N(x_s)} \lab{w_W} \ee
Moreover, for the canonical choice of the parameters $e_i$, i.e.
$e_3<e_2<e_1$ we have $u_n >0, \; n=1,2,\dots, N-1$ and hence all
the weights will be positive $w_s>0, \; s=0,1,\dots, N-1$. Finding
explicit expression for $x_s$ and $w_s$ is an interesting open
problem.

Note that we considered in this paper two special choices of the
parameters $q,\beta, \mu_1$ in expression \re{c0_FS} for $c_0(t)$
leading to positive discrete measures on the real axis. For
general choice of these parameters we obtain the function $c_0(t)$
which is quasi-periodic on the real axis $c_0(t+T) = e^{\nu}
c_0(t)$ where $T=2 \omega/w$ with some constant $\nu$. Using
Proposition {\bf 1} and its extension for quasi-periodic functions
(see the last paragraph of the first section) we can construct
orthogonal polynomials $P_n(x;t)$ with a measure located on some
horizontal line of the complex domain. In order to restore an
explicit expression for the measure in this generic case we need
explicit expression for the Fourier series of the quasi-periodic
functions $c_0(t)$ of type \re{c0_FS}. These Fourier series were
obtained in \cite{Tsu_Zhe}. A natural question arises: describe
all possible choices of the parameters $q, \beta, \mu_1$ when the
orthogonality measure is positive and located on the real axis.
This question demands a separate investigation.

We also mention interesting papers \cite{IVY}, \cite{IsMas} and
\cite{CF} where some explicit generalizations of the Stieltjes
continued fractions (and corresponding orthogonal polynomials)
were proposed. It is quite probable that these examples can be
considered as a special case (i.e. $t=0$) of more general elliptic
Toda chain solutions. It would be interesting to construct
corresponding explicit orthogonal polynomials and their
orthogonality measures. The author is indebted to M.Ismail for
bringing his attention to these papers.

\vspace{5mm}

{\Large \bf Acknowledgments.} \vspace{5mm}

The author thanks M.Ismail, A.Magnus, V.Spiridonov, S.Tsujimoto,
A.Veselov and L.Vinet for discussions of results of the papers.

\bb{99}

\bi{Akh_mom} N.I. Akhiezer [Achieser], The Classical Moment
Problem, Oliver and Boyd, Edinburgh, 1969 (originally published
Moscow, 1961).

\bi{Akhiezer} N.I. Akhiezer, {\it Elements of the Theory of
Elliptic Functions}, 2nd edition, "Nauka", Moscow, 1970.
Translations Math. Monographs {\bf 79}, AMS, Providence, 1990.

\bi{AMA} R. \'Alvarez-Nodarse, F. Marcell\'an, {\it Difference equation for modifications
of Meixner polynomials}.
J. Math. Anal. Appl. {\bf 194} (1995), 250--258.

\bi{Apt} A.I. Aptekarev, A. Branquinho, and F. Marcellan, {\it
Toda-type differential equations for the recurrence coefficients
of orthogonal polynomials and Freud transformation}. J. Comput.
Appl. Math. {\bf 78} (1997), no. 1, 139--160.

\bi{Askey} R.Askey, {\it Difference equation for modification of Meixner polynomials}.
In: Orthogonal Polynomials and Their Applications, (C.Brezinski et al. Eds). p. 418,
Annals of Computing and Applied Mathematics, Vol. 9, Baltzer AG Scientific, Basel, 1991.

\bi{Bavinck} H. Bavinck and H. Van Haeringen, {\it Difference equations for generalized Meixner polynomials},
J. Math. Anal. Appl., {\bf 184} (1994), 453--463.

\bi{Carlitz} L.Carlitz, {\it Some orthogonal polynomials related
to elliptic functions}. Duke Math. J. {\bf 27} (1960), 443–-459

\bi{Chi} T. Chihara, {\it An Introduction to Orthogonal Polynomials},
Gordon and Breach, NY, 1978.

\bi{ChiHamb} T.S.Chihara, {\it Hamburger moment problems and
orthogonal polynomials}, Transactions of the American Mathematical
Society, {\bf 315}, No. 1. (1989), 189--203.

\bi{Chu1} D.V. Chudnovsky and G.V. Chudnovsky, {\it Computer assisted number theory with applications}.
Number theory (New York, 1984--1985), 1--68, Lecture Notes in Math., {\bf 1240}, Springer, Berlin, 1987.

\bi{Chu2} D.V. Chudnovsky and G.V. Chudnovsky, {\it Transcendental methods and theta-functions},
Proc. Sympos. Pure Math. 49 (1989), 167-–232.

\bi{Chu}  D.V. Chudnovsky, G.V. Chudnovsky, {\it Hypergeometric
and modular function identities, and new rational approximations
to and continued fraction expansions of classical constants and
functions}. A tribute to Emil Grosswald: number theory and related
analysis, 117--162, Contemp. Math. {\bf 143}, Amer. Math. Soc.,
Providence, RI, 1993.

%\bi{Fr} F. G. Frobenius, {\it \"Uber die elliptischen Functionen
%zweiter Art}, J. Reine Angew. Math. 93 (1882), 53–-68.

%\bi{FrSt} F. G. Frobenius and L. Stickelberger, Zur Theorie der
%elliptischen Functionen, J. Reine Angew. Math. 83

\bi{CF} Eric van Fossen Conrad and Philippe Flajolet, {\it The Fermat cubic, elliptic functions, continued fractions, and a combinatorial excursion}, Seminaire Lotharingien de Combinatoire, (2006), {\bf 54} (B54g):1-44.

\bi{Ismail} M.E.H.Ismail, {\it Classical and Quantum orthogonal polynomials in one variable}.
Encyclopedia of Mathematics and its Applications (No. 98), Cambridge, 2005.

\bi{IVY} M. E. H. Ismail, G.  Valent and G. Yoon, {\it Some orthogonal polynomials related to elliptic functions},
J. Approx. Theory, {\bf 112} (2001), 251--178.

\bi{IsMas} M. E. H. Ismail and D. R. Masson, {\it Some continued fractions related to elliptic functions},  in ``Continued Fractions: From Analytic Theory to Constructive Approximation", B. C. Berndt and F. Gesztesy, eds.
Contemporary Mathematics,  {\bf 236}, 1999, pp. 149--166.

\bi{KS} Koekoek R and Swarttouw R F 1994 {\it The Askey scheme of
hypergeometric orthogonal \\polynomials and its q-analogue}, Report 94-05,
Faculty of Technical Mathematics and Informatics, Delft University of
technology.

\bi{Krall} A.M.Krall, {\it Hilbert space, boundary value problems and orthogonal polynomials}.
Operator Theory: Advances and Applications, {\bf 133}. Birkhauser Verlag, Basel, 2002

\bi{LoBr} J.S.Lomont, J. Brillhart, {\it Elliptic polynomials}.
Chapman \& Hall/CRC, Boca Raton, FL, 2001.

\bi{Lukacs} E. Lukacs, {\it Characteristic functions}, 2ed.,
Griffin, 1970.

\bi{Milne} S.Milne, {\it Infinite Families of Exact Sums of
Squares Formulas, Jacobi Elliptic Functions, Continued Fractions,
and Schur Functions}, Ramanujan J., {\bf 6} (2002), 7-–149,

\bi{NaZhe} Y.Nakamura and A.Zhedanov, {\it Special solutions of
the Toda chain and combinatorial numbers}, J. Phys. A: Math. Gen.
{\bf 37}, (2004), 5849-–5862.

\bi{Peh1} F. Peherstorfer, {\it On Toda lattices and orthogonal
polynomials}.
Proceedings of the Fifth International Symposium on Orthogonal
Polynomials,
Special Functions and their Applications (Patras, 1999). J. Comput.
Appl. Math. {\bf 133}
  (2001), 519--534.

\bi{Peh2} F. Peherstorfer, V. Spiridonov and A. Zhedanov, {\it
The Toda chain, the Stieltjes function, and orthogonal polynomials}.
(Russian) Teoret. Mat. Fiz. {\bf 151} (2007), no. 1, 81--108.

\bi{Rogers}  L.J.Rogers, {\it On the Representation of Certain
Asymptotic Series as Convergent Continued Fractions},  Proc.
London Math. Soc. (1907) s2-4: 72--89.

\bi{ShoTa} J. Shohat and J. D. Tamarkin, {\it The problem of
moments,} Math. Surveys, no. 1, Amer. Math. Soc., Providence,
R.I., 1943/1950.

\bi{Sogo} K. Sogo, {\it Time-dependent orthogonal polynomials and
theory of soliton. Applications to matrix model, vertex model and
level statistics}. J. Phys. Soc. Japan {\bf 62} (1993), 1887--1894

%\bi{Sz} G. Szeg\H{o}, Orthogonal Polynomials, fourth edition,  AMS,
%1975.

\bibitem{Toda} M. Toda, {\it Theory of Nonlinear Lattices}.
Second edition. Springer Series in Solid-State Sciences, {\bf 20},
Springer-Verlag, Berlin, 1989. x+225 pp.

\bi{Tsu_Zhe} S.Tsujimoto and A.Zhedanov, {\it Elliptic
hypergeometric Laurent biorthogonal polynomials with a dense point
spectrum on the unit circle}, SIGMA {\bf 5} (2009), 033, arXiv:0809.2574,

\bi{Valent3} G.Valent, {\it Asymptotic analysis of some associated
orthogonal polynomials connected with elliptic functions}, SIAM J.
Math.Anal., {\bf 25} (1994), 749--775.

%\bi{Valent2} {\it  Exact solutions of processes some quadratic and
%quartic birth and death and related orthogonal polynomials},  J.
%Comput. Appl. Math. {\bf 67} (1996), 103--127.

\bi{Valent1} G.Valent, {\it Associated Stieltjes-Carlitz polynomials and a
generalization of Heun's differential equation}, J. Comput. Appl.
Math. {\bf 57} (1995), 293--307.

\bi{Valent} G.Valent, {\it From asymptotics to spectral measures:
determinate versus indeterminate moment problems}, Mediterr. J.
Math. {\bf 3} (2006), 327-–345.

\bi{WW} E.T. Whittacker, G.N. Watson, {\em A Course of Modern
Analysis}, Cambridge, 1927.

%\bi{Wimp}  J. Wimp, {\it Hankel determinants of some polynomials
%arising in combinatorial analysis}. Computational methods from
%rational approximation theory (Wilrijk, 1999). Numer. Algorithms
%{\bf 24} (2000), no. 1-2, 179--193.

\bi{Zhe_cndn} A. Zhedanov, {\it Elliptic polynomials orthogonal on
the unit circle with a dense point spectrum}, Ramanujan J. (2009) {\bf 19}, 351-–384, arXiv:0711.4696

\eb

\end{document}